\documentclass[12pt]{amsart}
\usepackage{amsmath,amsfonts,amsthm,amssymb,color}

\usepackage[T1]{fontenc}
\usepackage{mathrsfs}
\usepackage{pdfsync}

\newcommand{\bd}{\mathbf{D}}
\newcommand{\ou}{[0,1]}

\newcommand{\bb}{\mathbf{B}}

\newcommand{\der}{\delta}

\newcommand{\ha}{\hat{a}}
\newcommand{\hr}{\hat{r}}
\newcommand{\hz}{\hat{z}}
\newcommand{\hze}{\hat{\zeta}}
\newcommand{\id}{\mbox{Id}}

\newcommand{\ist}{\int_{s}^{t}}
\newcommand{\norm}[1]{\lVert #1\rVert}

\newcommand{\ott}{[0,T]}
\newcommand{\xd}{\mathbf{B}^{\mathbf{2}}}
\newcommand{\1}{{\bf 1}}
\newcommand{\2}{{\bf 2}}

\newcommand{\lp}{\left(}
\newcommand{\rp}{\right)}
\newcommand{\lc}{\left[}
\newcommand{\rc}{\right]}
\newcommand{\lcl}{\left\{}
\newcommand{\rcl}{\right\}}
\newcommand{\lln}{\left|}
\newcommand{\rrn}{\right|}

\newcommand{\ep}{\varepsilon}

\newcommand{\ga}{\gamma}
\newcommand{\ka}{\kappa}
\newcommand{\la}{\lambda}
\newcommand{\laa}{\Lambda}
\newcommand{\om}{\omega}
\newcommand{\oom}{\Omega}

\newcommand{\ze}{\zeta}
\newcommand{\vp}{\varphi}

\newcommand{\beq}{\begin{equation}}
\newcommand{\eeq}{\end{equation}}
\newcommand{\bea}{\begin{eqnarray}}
\newcommand{\eea}{\end{eqnarray}}
\newcommand{\beas}{\begin{eqnarray*}}
\newcommand{\eeas}{\end{eqnarray*}}

\def\msf{{\mathscr F}}

\def\me{{\mathbb  E}}

\def\mr{{\mathbb  R}}

\def\mp{{\mathbb  P}}

\newcommand{\cb}{{\mathcal B}}
\newcommand{\cac}{{\mathcal C}}

\newcommand{\cf}{{\mathcal F}}

\newcommand{\ch}{{\mathcal H}}
\newcommand{\cj}{{\mathcal J}}

\newcommand{\cm}{{\mathcal M}}
\newcommand{\cn}{{\mathcal N}}

\newcommand{\cq}{{\mathcal Q}}

\newcommand{\cs}{{\mathcal S}}

\newcommand{\cz}{{\mathcal Z}}

\newcommand{\D}{{\mathbb D}}
\newcommand{\EE}{{\mathbb E}}

\newcommand{\PP}{{\mathbb P}}

\newcommand{\R}{{\mathbb R}}

\newtheorem{theorem}{Theorem}[section]

\newtheorem{corollary}[theorem]{Corollary}

\newtheorem{definition}[theorem]{Definition}

\newtheorem{hypothesis}[theorem]{Hypothesis}
\newtheorem{lemma}[theorem]{Lemma}

\newtheorem{proposition}[theorem]{Proposition}

\theoremstyle{remark}
\newtheorem{remark}[theorem]{Remark}
\theoremstyle{remark}
\newtheorem{example}[theorem]{Example}

\newtheorem{foo}[theorem]{Remarks}

\title[Gaussian bounds for fractional SDEs]{ Upper bounds for the density of  solutions of stochastic differential equations driven by fractional Brownian motions}

\author[F. Baudoin \and C. Ouyang \and S. Tindel]{Fabrice Baudoin \and Cheng Ouyang \and Samy Tindel}

\thanks{First author supported in part by
NSF Grant DMS 0907326. Third author partially supported by the (French) ANR grant ECRU}

\address{Fabrice Baudoin, Dept. Mathematics, Purdue University, 150 N. University St., West Lafayette, IN 47907-2067, USA.}
\email{fbaudoin@math.purdue.edu}

\address{Cheng Ouyang, Dept. Mathematics, Purdue University, 150 N. University St., West Lafayette, IN 47907-2067, USA.}
\email{couyangmath@gmail.com}

\address{Samy Tindel, Institut {\'E}lie Cartan Nancy, Universit\'e de Nancy 1, B.P. 239,
54506 Vand{\oe}uvre-l{\`e}s-Nancy Cedex, France.}
\email{tindel@iecn.u-nancy.fr}

\begin{document}

\maketitle

\

\begin{abstract}
In this paper we study upper bounds for the density of solution of stochastic differential equations driven by a fractional Brownian motion with Hurst parameter $H> 1/3$. We show that under some geometric conditions, in the regular case $H>1/2$, the density of the solution satisfy the log-Sobolev inequality, the Gaussian concentration inequality and  admits an upper Gaussian bound. In the rough case $H>1/3$ and under the same geometric conditions, we show that the density of the solution is smooth and admits an upper sub-Gaussian bound.
\end{abstract}

\tableofcontents

\section{Introduction}

Let $B=(B^1,\ldots,B^d)$ be a $d$ dimensional fractional Brownian motion (fBm in the sequel) defined on a complete probability space $(\Omega,\cf,\PP)$, with Hurst parameter $H\in(0,1)$. Recall that it means that $B$ is a centered Gaussian process indexed by $\R_+$, whose coordinates are independent and satisfy
\begin{equation}\label{eq:var-increm-fbm}
\EE\lc  \lp B_t^j -B_s^j\rp^2\rc= |t-s|^{2H}, \quad\mbox{for}\quad s,t\in\R_+.
\end{equation}
In particular, by considering the family $\{B^H;\, H\in (0,1)\}$, one obtains some Gaussian processes with any prescribed Hölder regularity, while fulfilling some intuitive scaling properties. This converts fBm into the most natural generalization of Brownian motion to this day.

We are concerned here with the following class of equations driven by $B$:
\begin{equation}
\label{eq:sde-intro} X^{x}_t =x +\int_0^t V_0 (X^x_s)ds+
\sum_{i=1}^d \int_0^t V_i (X^{x}_s) dB^i_s,
\end{equation}
where $x$ is a generic initial condition and $\{V_i;\, 0\le i\le d\}$ is a collection of smooth vector fields of $\R^d$. Owing to the fact that fBm is a natural generalization of Brownian motion, this kind of model is often used by practitioners in different contexts, among which we would like to highlight recent sophisticated models  in Biophysics~\cite{KS,SW,TBV}.

As far as mathematical results are concerned, equation (\ref{eq:sde-intro}) is now a fairly well understood object: existence and uniqueness results are obtained for $H>\frac{1}{2}$ thanks to Young integral type tools \cite{Za,NR}, while rough paths methods \cite{FV-bk,LQ} are required for $\frac{1}{4}<H<\frac{1}{2}$. Numerical schemes can be implemented for this kind of systems \cite{DNT,FV-bk}, and a notion of ergodicity is also available \cite{Ha,HO}. Finally, the law of $X_t^x$ has been analyzed by means of semi-group type methods \cite{BC,NNRT} and its density has also been investigated in \cite{BH,CFV,HN,NS}. 

In spite of these advances, concentrations results and Gaussian bounds for the solution to (\ref{eq:sde-intro}) are scarce: we are only aware of the large deviation results \cite{MS} in this line of investigation. The current article is thus an attempt to make a step in this direction, by analyzing a special but nontrivial situation.

Indeed, we consider here equation (\ref{eq:sde-intro}) driven by a fBm with Hurst parameter $H\in(\frac{1}{3}, 1)$, and we suppose that our vector fields $V_0,\ldots,V_d$ fulfill either of the following non-degeneracy and antisymmetric hypothesis:
\begin{hypothesis}\label{H}
The vector fields $V_0,\ldots,V_d$ are $C^\infty$-bounded, and $V_1,\ldots,V_d$ satisfy

\smallskip
\noindent
\emph{(i)}
For every $x \in \mathbb{R}^d$, the vectors $V_1(x),\cdots,V_d(x)$ form a basis of $\mathbb{R}^d$.

\smallskip

\noindent
\emph{(ii)}
There exist smooth and bounded functions $\omega_{ij}^k$ such that:
\begin{equation}\label{eq:antisym-V}
[V_i,V_j]=\sum_{k=1}^d \omega_{ij}^k V_k,
\quad\mbox{and}\quad
\omega_{ij}^k =-\omega_{ik}^j.
\end{equation}
\end{hypothesis}

The second assumption (ii) is of geometric nature and actually means that the Levi-Civita connection associated with the Riemannian structure given by the vector fields $V_i$'s is 
\[
\nabla_X Y=\frac{1}{2} [X,Y].
\]
In a Lie group structure, this is equivalent to the fact that the Lie algebra is of compact type, or in other words that the adjoint representation is unitary. Such geometric assumption already appeared in the work  \cite{BO} where it was used to prove  a small-time asymptotics of the density.

\begin{hypothesis}\label{H'}
The Hypothesis \ref{H} is satisfied and moreover, the vector fields $V_1,...,V_d$ form moreover a uniform elliptic system. That is
$$|v^T VV^Tv|\geq \lambda|v|^2,\quad\quad\mathrm{for\ all}\ v\in \mr^d.$$
Here  $V=(V^i_j)_{i,j=1,...,d}$ and $\lambda$ is a positive constant.
\end{hypothesis}

When $H>\frac{1}{2}$, under Hypothesis \ref{H} our main result can be loosely summarized as follows (see Theorem \ref{thm:main-gauss-bnd} for a precise statement):
\begin{theorem}\label{thm:main-intro}
Fix $H>\frac{1}{2}$. Let $X^x$ be the solution to equation (\ref{eq:sde-intro}), and suppose Assumption \ref{H} is satisfied. Then for any $t\in\R_+^*$, the random variable $X_t^x$ admits a smooth density $p_X(t,\cdot)$. Furthermore, there exist 3 positive constants $c_t^{(1)},c_t^{(2)},c_{t,x}^{(3)}$ such that
\begin{equation*}
p_X(t,y) \le c_t^{(1)} \exp\lp - c_t^{(3)} \lp  |y|- c_{t,x}^{(2)} \rp^2\rp,
\end{equation*}
for any $y\in\R^d$.
\end{theorem}
We don't claim any optimality in the quantities $c_t^{(1)},c_t^{(2)}$ and $c_{t,x}^{(3)}$ above (whose exact definitions are postponed to Section \ref{sec:gauss-bnd}). Nevertheless, this is (to the best of our knowledge) the first Gaussian type bound available for solutions of differential equations driven by fBm.

Let us say a few words about the strategy we have followed in order to prove Theorem~\ref{thm:main-intro}. It is mostly based on stochastic analysis tools, and particularly on a general integration by parts formula giving an exact expression for the density $p_X(t,\cdot)$ in terms of Malliavin derivatives in the non-degenerate case we are dealing with. In this context, it is crucial to bound the first Malliavin derivative of $X_t^x$ (called $\mathbf{D}X_t^x$ in the sequel) efficiently. This is where our asymmetry hypothesis on the vector fields $V_1,\ldots,V_n$ enter into the picture, and we shall see (at Theorem \ref{bound MalliavinD})  how asymmetry properties yield an easy deterministic bound on $\bd X_t^x$. This result enables to get concentration results for the law of $X_t^x$, and is the key to our density bounds as well.

As another interesting consequences of the deterministic bound on $\mathbf{D}X^x_t$, we also obtain Log-Sobolev inequality and Poincar\'{e} inequality for the law of $X_t$.

\

Once the picture for the smooth case (when $H>\frac{1}{2}$) becomes clear, we are able to extend some of our results described above to the irregular case when $\frac{1}{3}<H<\frac{1}{2}$. In particular, we are able to prove 

\begin{theorem}
Fix $H\in (\frac{1}{3}, \frac{1}{2})$. Assume Hypothesis \ref{H'}.  Let $X^x$ be the solution to equation (\ref{eq:sde-intro}) and $\gamma_{X_t}$ the Malliavin matrix of $X^x_t$, $t>0$. We have 
$|\det\gamma_{X_t}|^{-1}\in L^{\infty}(\mp)$. The random variable $X_t^x$ admits a smooth density $p_X(t,\cdot)$ and for any $\delta<H$ there exist 2 positive constants $c_t^{(1)},c_t^{(2)}$ such that
\begin{equation}\label{eq:sub-gauss-bnd}
p_X(t,y) \le c_t^{(1)} \exp\lp - c_t^{(2)}  |y|^{\delta} \rp,
\end{equation}
for all $y\in\R^d$.
\end{theorem}
The existence of a density for solutions to stochastic differential equations of the form (\ref{eq:sde-intro}) under H\"{o}rmander's condition has been obtained by Cass and Friz \cite{CF} for any $\frac{1}{4}<H<\frac{1}{2}$.  While finishing the current article, an important step towards the study of regular densities in the rough case $\frac{1}{4}<H<\frac{1}{2}$ has been accomplished in \cite{CLL}, where integrability estimates for the Jacobian of equation (\ref{eq:sde-intro}) are established.
Nevertheless, as of today, besides the result of P. Driscoll \cite{Dr}, when $H<\frac{1}{2}$, to the best of our knowledge, Hypothesis \ref{H'} is a first wide class of examples where we have an affirmative answer for the smoothness of the density. Our bound on the inverse of the Malliavin matrix together with polynomial bounds on the H\"older norm of the Malliavin derivative allows then to obtain the sub-Gaussian upper bound (\ref{eq:sub-gauss-bnd}).


\smallskip

\noindent
\textbf{Notations:} Throughout this paper, unless otherwise specified we use $|\cdot|$ for Euclidean norms and $\|\cdot\|_{L^p}$ for the $L^p$ norm with respect to the underlying probability measure $\mp$.

Consider a finite-dimensional vector space $V$. The space of $V$-valued Hölder continuous functions defined on $[0,1]$, with Hölder continuity exponent $\ga\in(0,1)$, will be denoted by $\cac^\ga(V)$, or just $\cac^\ga$ when this does not yield any ambiguity. For a function $g\in\cac^\ga(V)$ and $0\le s<t\le 1$, we shall consider the semi-norms
\begin{equation}\label{eq:def-holder-norms}
\|g\|_{s,t,\ga}=\sup_{s\le u<v\le t}\frac{|g_v-g_u|_{V}}{|v-u|^{\ga}},
\end{equation}
The semi-norm $\|g\|_{0,1,\ga}$ will simply be denoted by $\|g\|_{\ga}$.

\section{Stochastic calculus for fractional Brownian motion}
\label{sec:stoch-calc-fbm}
For some fixed $H\in(\frac{1}{3},1)$, we consider $(\oom,\cf,\PP)$ the canonical probability space associated with the fractional
Brownian motion (in short fBm) with Hurst parameter $H$. That is,  $\oom=\cac_0([0,1])$ is the Banach space of continuous functions
vanishing at $0$ equipped with the supremum norm, $\cf$ is the Borel sigma-algebra and $\PP$ is the unique probability
measure on $\oom$ such that the canonical process $B=\{B_t=(B^1_t,\ldots,B^d_t), \; t\in [0,1]\}$ is a fractional Brownian motion with Hurst
parameter $H$.
In this context, let us recall that $B$ is a $d$-dimensional centered Gaussian process, whose covariance structure is induced by equation~(\ref{eq:var-increm-fbm}). This can be equivalently stated as 
\begin{equation*}
R\left( t,s\right) :=\EE\lc  B_s^j \, B_t^j\rc
=\frac{1}{2}\left( s^{2H}+t^{2H}-|t-s|^{2H}\right),
\quad\mbox{for}\quad 
s,t\in[0,1] \mbox{ and } j=1,\ldots,d.
\end{equation*}
In particular it can be shown, by a standard application of Kolmogorov's criterion, that $B$ admits a continuous version
whose paths are $\ga$-H\"older continuous for any $\ga<H$. 

This section is devoted to give the basic elements of stochastic calculus with respect to $B$ which allow to understand the remainder of the paper.

\subsection{Malliavin calculus tools}\label{sec:malliavin-tools}
Gaussian techniques are obviously essential in the analysis of fBm, and we proceed here to introduce some of them (see \cite{Nu06} for further details): let $\mathcal{E}$ be the space of $\mathbb{R}^d$-valued step
functions on $[0,1]$, and $\mathcal{H}$  the closure of
$\mathcal{E}$ for the scalar product:
\[
\langle (\mathbf{1}_{[0,t_1]} , \cdots ,
\mathbf{1}_{[0,t_d]}),(\mathbf{1}_{[0,s_1]} , \cdots ,
\mathbf{1}_{[0,s_d]}) \rangle_{\mathcal{H}}=\sum_{i=1}^d
R(t_i,s_i).
\]
We denote by $K^*_H$ the isometry between $\mathcal{H}$ and $L^2([0,1])$.
When $H>\frac{1}{2}$ it can be shown that $\mathbf{L}^{1/H} ([0,1], \mathbb{R}^d)
\subset \mathcal{H}$, and when $\frac{1}{3}<H<\frac{1}{2}$ one has 
$$C^\gamma\subset \mathcal{H}\subset L^2([0,1])$$
for all $\gamma>\frac{1}{2}-H$.

Some isometry arguments allow to define the Wiener integral $B(h)=\int_0^{1} \langle h_s, dB_s \rangle$ for any element $h\in\ch$, with the additional property $\EE[B(h_1)B(h_2)]=\langle h_1,\, h_2\rangle_{\ch}$ for any $h_1,h_2\in\ch$.
A $\mathcal{F}$-measurable real
valued random variable $F$ is then said to be cylindrical if it can be
written, for a given $n\ge 1$, as
\begin{equation*}
F=f\lp  B(h^1),\ldots,B(h^n)\rp=
f \Bigl( \int_0^{1} \langle h^1_s, dB_s \rangle ,\ldots,\int_0^{1}
\langle h^n_s, dB_s \rangle \Bigr)\;,
\end{equation*}
where $h^i \in \mathcal{H}$ and $f:\mathbb{R}^n \rightarrow
\mathbb{R}$ is a $C^{\infty}$ bounded function with bounded derivatives. The set of
cylindrical random variables is denoted $\mathcal{S}$. 

The Malliavin derivative is defined as follows: for $F \in \mathcal{S}$, the derivative of $F$ is the $\mathbb{R}^d$ valued
stochastic process $(\mathbf{D}_t F )_{0 \leq t \leq 1}$ given by
\[
\mathbf{D}_t F=\sum_{i=1}^{n} h^i (t) \frac{\partial f}{\partial
x_i} \left( B(h^1),\ldots,B(h^n)  \right).
\]
More generally, we can introduce iterated derivatives. If $F \in
\mathcal{S}$, we set
\[
\mathbf{D}^k_{t_1,\ldots,t_k} F = \mathbf{D}_{t_1}
\ldots\mathbf{D}_{t_k} F.
\]
For any $p \geq 1$, it can be checked that the operator $\mathbf{D}^k$ is closable from
$\mathcal{S}$ into $\mathbf{L}^p(\oom;\mathcal{H}^{\otimes k})$. We denote by
$\mathbb{D}^{k,p}(\mathcal{H})$ the closure of the class of
cylindrical random variables with respect to the norm
\[
\left\| F\right\| _{k,p}=\left( \mathbb{E}\left( F^{p}\right)
+\sum_{j=1}^k \mathbb{E}\left( \left\| \mathbf{D}^j F\right\|
_{\mathcal{H}^{\otimes j}}^{p}\right) \right) ^{\frac{1}{p}},
\]
and
\[
\mathbb{D}^{\infty}(\mathcal{H})=\bigcap_{p \geq 1} \bigcap_{k
\geq 1} \mathbb{D}^{k,p}(\mathcal{H}).
\]

\subsection{Differential equations driven by fBm}
Recall that we consider the following kind of equation:
\begin{equation}
\label{eq:sde} X^{x}_t =x +\int_0^t V_0 (X^x_s)ds+
\sum_{i=1}^d \int_0^t V_i (X^{x}_s) dB^i_s,
\end{equation}
where the vector fields $V_0,\ldots,V_n$ are $C^\infty$-bounded.

When equation (\ref{eq:sde}) is driven by a fBm with Hurst parameter $H>\frac{1}{2}$ it can be solved, thanks to a fixed point argument, with the stochastic integral interpreted in the  (pathwise) Young sense (see e.g. \cite{Gu}). Let us recall that Young's integral can be defined in the following way:
\begin{proposition}\label{prop:young-intg}
Let $f\in \cac^\ga$, $g\in \cac^\kappa$ with $\ga+\kappa>1$, and $0\le s\le t\le 1$. Then the  integral $\int_s^t g_\xi df_\xi$ is well-defined as limit of Riemann sums along partitions of $[s,t]$. Moreover, the following estimation is fulfilled:
\beq
\left|\int_s^t g_\xi df_\xi\right| \leq C \|f\|_\ga \|g\|_\kappa |t-s|^\ga,
\label{y}
\eeq
where the constant $C$ only depends on $\ga$ and $\kappa$. A sharper estimate is also available:
\begin{equation}\label{eq:ineq-young-sharp}
\left|\int_s^t g_\xi df_\xi\right| \leq |g_s| \, \|f\|_\ga |t-s|^\ga
+ c_{\ga,\ka}  \|f\|_\ga \|g\|_\kappa |t-s|^{\ga+\ka}.
\end{equation}
\end{proposition}

\smallskip

With this definition in mind, we can solve our differential system of interest, and the following moments bounds are proven in \cite{HN}:
\begin{proposition}\label{prop:moments-sdes}(Hu-Nualart)
Consider equation (\ref{eq:sde}) driven by a fBm $B$ with Hurst parameter $H>\frac{1}{2}$. Let us call $X^x$ its unique $\beta$-H\"older continuous solution, for any $\beta<H$. Then

\smallskip

\noindent
(1) When vector fields $V$ are $C^\infty$-bounded, we have
\begin{equation*}
\sup_{t\in[0,T]} |X_t^x| \le |x|+ c_{V,T} \|B\|_{0,T,\beta}^{1/\beta}.
\end{equation*}

\smallskip

\noindent
(2) If we only assume that vector fields $V$ have linear growth, with $\nabla V, \nabla^2 V$, bounded, the following estimate holds true:
\begin{equation}\label{eq:bnd-lin-sde}
\sup_{t\in[0,T]} |X_t^x| \le \lp 1+ |x|\rp \,  \exp\lp c_{V,T} \|B\|_{0,T,\beta}^{1/\beta}\rp.
\end{equation}
\end{proposition}

\begin{remark}
The framework of fractional integrals is used in \cite{HN} in order to define integrals with respect to $B$. It is however easily seen to be equivalent to the Young setting we have chosen to work with.
\end{remark}

When the Hurst parameter $\frac{1}{3}<H<\frac{1}{2}$, equation (\ref{eq:sde}) can be solved, again by fixed point argument, with the stochastic integral interpreted in the  (pathwise) rough path theory (see e.g. \cite{Gu} and \cite{LQ}).
In this case, we obtain
\begin{proposition}\label{prop:moments-sdes-rough}(Besal\'u-Nualart \cite{BN})
Consider equation (\ref{eq:sde}) driven by a fBm $B$ with Hurst parameter $\frac{1}{3}<H<\frac{1}{2}$. Denote by $X^x$ its unique $\beta$-H\"older continuous solution, for any $\beta<H$. 
If the vector fields $V$ are $C^\infty$-bounded, then for any $\lambda>0$ and $\delta<H$
\begin{equation*}
\me\left(\exp\lambda\left(\sup_{0\leq t\leq T}|X_t|^\delta\right)\right)<\infty.
\end{equation*}
\end{proposition}

Once equation (\ref{eq:sde}) is solved, the vector $X_t^x$ is a typical example of  random variable which can be differentiated in the Malliavin sense. In fact, fix $H\in(\frac{1}{3},1)$, one gets the following results (see \cite{CF} and \cite{NS}  for further details):
\begin{proposition}\label{prop:deriv-sde}
Let $X^x$ be the solution to equation (\ref{eq:sde}) and suppose $V_i$'s are $C^\infty$-bounded vector fields on $\mr^d$. Then
for every $i=1,\ldots,d$, $t>0$, and $x \in \mathbb{R}^d$, we have $X_t^{x,i} \in
\mathbb{D}^{\infty}(\mathcal{H})$ and
\begin{equation*}
\mathbf{D}^j_s X_t^{x}= \mathbf{J}_{0 \rightarrow t}
\mathbf{J}_{0 \rightarrow s}^{-1} V_j (X_s) ,~~j=1,\ldots,d, ~~ 0\leq
s \leq t,
\end{equation*}
where $\mathbf{D}^j_s X^{x,i}_t $ is the $j$-th component of
$\mathbf{D}_s X^{x,i}_t$, and $\mathbf{J}_{0\rightarrow t}=\frac{\partial X^x_t}{\partial x}$. 

\end{proposition}

Finally the following approximation result, which can be found for instance in \cite{FV-bk}, will also be used in the sequel:
\begin{proposition}\label{prop:lin-interpol}
For $m\ge 1$ and $T>0$, let
$B^{m}=\{B^m_t;\, t\in [0,T]\}$ be the sequence of linear interpolations of
$B$ along the dyadic subdivision of $[0,T]$ of mesh $m$; that
is if $t_i^m= i 2^{-m} T$ for $i=0,..., 2^m;$ then for $t \in
(t_i^m, t_{i+1}^m ],$
\begin{equation*}
B^{m}_t=B_{t_{i}^{m}} +
\frac{t-t_{i^m}}{t_{i+1}^m - t_i^m} (B_{t_{i+1}^m}
-B_{t_{i}^m} ).
\end{equation*}
Consider $X^{m}$ the solution to equation (\ref{eq:sde}) restricted to $[0,T]$, where $B$ has been replaced by $B^{m}$. Set also $\mathbf{D}^j_s X_t^{m}= \mathbf{J}_{0 \rightarrow t}\mathbf{J}_{0 \rightarrow s}^{-1} V_j (X_s^{m})$, for $j=1,\ldots,d$ and $0\leq s \leq t$. Then almost surely, for any $\ga<H$ and $t\in[0,T]$ the following holds true:
\begin{equation}\label{eq:limit-interpol}
\lim_{m\to\infty} \lp \|X^{x}-X^{m}\|_{\ga} + \|\mathbf{D}^j X_t^{x}-\mathbf{D}^j X_t^{m}\|_{\ga} \rp=0.
\end{equation}
\end{proposition}


\section{ Estimates  for solutions of SDEs driven by fBm: the smooth case}

Throughout this section, we fix $H\in (\frac{1}{2}, 1)$. Recall that $X^x$ designates the solution to (\ref{eq:sde}). This section is devoted to get some further bounds for $X_t^x$ and its Malliavin derivatives, under Assumption \ref{H}. 

Notice that among our set of hypothesis, the antisymmetric property (\ref{eq:antisym-V}) for the vector fields $V_1,\ldots,V_n$ is the most specific one. It will be mainly used through the following lemma:
 \begin{lemma}\label{contraction}
Let $A_1,A_2$ be $n\times n$ matrices, whose exponential are defined by $e^{A_j}=\sum_{n=0}^{\infty}A_j^n/n!$ for $j=1,2$. If we assume that $A_2$ is skew symmetric, then
\[
\left\| e^{A_1+A_2} \right\| \le e^{\| A_1 \|},
\]
where $\|A\|$ stands for the Euclidean norm of a matrix $A$.
\end{lemma}

 \begin{proof}
 Let us first prove the following (presumably classical) identity:
 \begin{equation}\label{eq:duhamel}
 e^{t(A_1+A_2)}=e^{tA_2}-\int_0^t e^{(t-s)(A_1+A_2)} Ae^{sA_2} ds.
\end{equation}
Indeed, consider the function $s\mapsto\vp(s)$ defined on $[0,1]$ by $\vp(s)=e^{(t-s)(A_1+A_2)}e^{sA_2}$. Then it is easily seen that $\vp$ is differentiable and
\begin{equation*}
\vp'(s)= e^{(t-s)(A_1+A_2)} A_1e^{sA_2}.
\end{equation*}
 By writing
 \begin{equation*}
e^{tA_2}-e^{t(A_1+A_2)}=\vp(t)-\vp(0)=\int_0^t \vp'(s) \, ds,
\end{equation*}
relation (\ref{eq:duhamel}) is now easily obtained.

Let us see now the implications of  (\ref{eq:duhamel}): according to the fact that $A_2$ is skew-symmetric, we have $\|e^{sA_2}\|\le 1$ for any $s\ge 0$. Therefore,
 \[
 \left\| e^{t(A_1+A_2)} \right\| \le 1+\| A_1 \| \int_0^t  \left\|  e^{(t-s)(A_1+A_2)} \right\|  ds.
 \]
 By denoting $f(t)=\| e^{t(A_1+A_2)}\|$ we thus get 
 \[
 f(t) \le 1 +\| A_1 \| \int_0^t f(s) ds.
 \]
 This implies $f(t)\le e^{\| A_1 \| t}$ by a standard application of Gronwall's lemma and finishes our proof.

 \end{proof}

We are now ready to prove the main result of this section, which is an almost sure deterministic bound for the Malliavin derivative of $X_t^x$:

\begin{theorem}\label{bound MalliavinD}
Under Assumption \ref{H}, the Malliavin derivative of the solution $X_T^x$ to equation (\ref{eq:sde}) can be bounded as follows for any $T\in[0,1]$ (in the almost sure sense):
\begin{equation}\label{eq:bnd-first-deriv-infty}
\left\| \mathbf{D} X_T^{x} \right\|_{\infty}    \le M   \exp \left( CT \right),
\quad\mbox{with}\quad
M=\sup_{x \in \mathbb{R}^n} \sup_{\| \lambda \| \le 1}\left| \sum_{i=1}^d \lambda_i V_i(x) \right|^2,
\end{equation}
and where the constant $C$ linearly depends on $V_0$. In particular, one also has $\left\| \mathbf{D} X_T^{x} \right\|_{\ch}    \le M   \exp \left( CT \right)$.
\end{theorem}

\begin{proof}
Let us focus on the proof of (\ref{eq:bnd-first-deriv-infty}). Indeed, since $\|f\|_{\ch}$ is dominated by the supremum norm when $H>\frac{1}{2}$, this will be sufficient in order to prove the second claim of our theorem. We now split our proof in two steps.

\smallskip

\noindent
\textit{Step1: Matricial expression for the derivative.}
Let us first restate Proposition \ref{prop:deriv-sde} in the following form: $\mathbf{D} X_T^{x}$ is solution to 
\[
\mathbf{D}^j_s X_T^{x}=\mathbf{J}_{0 \rightarrow T} (\Phi^*_s V_j) (x), \quad 0\le s\le T,
\]
where $\Phi^*_s V_j$ denotes the pullback action of the diffeomorphism $\Phi_s=X_s^\cdot: \mr^d\rightarrow \mr^d$ on the vector field $V_j$. Now, a simple application of the change of variable formula for Young type integrals yields
\[
d (\Phi^*_s V_j) (x)=(\Phi^*_s [V_0,V_j]) (x)ds +\sum_{i=1}^d (\Phi^*_s [V_i,V_j]) (x)dB^i_s.
\]
Moreover, recall that the Lie brackets $[V_i,V_j]$ can be decomposed, according to Assumption~\ref{H}, into
\[
[V_0,V_j]=\sum_{k=1}^d \omega_{0j}^k V_k,
\quad\mbox{and}\quad
[V_i,V_j]=\sum_{k=1}^d \omega_{ij}^k V_k,
\]
with $\omega_{ij}^k=-\omega_{ik}^j$ for $i,j,k\ge 1$. Hence,
\[
d (\Phi^*_s V_j) (x)=\sum_{k=1}^d \omega_{0j}^k (X_s^x) (\Phi^*_s V_k) (x)ds +\sum_{k=1}^d\sum_{i=1}^d \omega_{ij}^k (X_s^x) (\Phi^*_s V_k) (x)dB^i_s.
\]
By denoting  $ \mathcal{M}_s$ the $d\times d$ matrix with columns 
\[
\mathcal{M}^j_s= \mathbf{D}^j_s X_T^{x}=\mathbf{J}_{0 \rightarrow T} (\Phi^*_s V_j) (x),
\]
we therefore obtain the equation
\begin{align} \label{eq:ed-def}
d\mathcal{M}_s=\mathcal{M}_s \left( \omega_{0} (X_s^x) ds + \sum_{i=1}^d \omega_{i} (X_s^x) dB^i_s\right), \quad \mathcal{M}_T=V(X^x_T),
\end{align}
where $V(X_T^x)$ is the matrix with columns $V_j (X^x_T)$, $1 \le j \le d$ and where  $ \omega_{i} (X_s^x)$ is the skew symmetric matrix with entries  $\omega_{ij}^k (X_s^x)$.

\smallskip

\noindent
\textit{Step2: Approximation procedure.}
In order to show that  the process $\mathcal{M}$ is uniformly bounded, consider the dyadic approximation introduced at Proposition \ref{prop:lin-interpol}. By applying (\ref{eq:limit-interpol}) to the couple $(X, \cm)$, it is sufficient to prove our uniform bounds on $\mathbf{D}^j X_T^{m}$, uniformly in $m$. Let us thus consider  $\mathcal{M}^m$ the solution of (\ref{eq:ed-def}) where $B$ is
replaced by $B^{m}$, that is
$$
d\mathcal{M}^m_s=\mathcal{M}^m_s \left( \omega_{0} (X_s^m) ds + \sum_{k=1}^d \omega_{k} (X_s^m) dB^{m,k}_s\right), \quad \mathcal{M}^m_T=V(X^m_T).
$$
In the sequel set also 
\begin{equation*}
\Delta B_{{t_{n-1}^{m}}t_{n}^{m}}^{k}:=
\frac{B^k_{t_{n}^{m}}-B^k_{t_{n-1}^{m}} }{t_{n}^m - t_{n-1}^m},
\quad\mbox{for}\quad 1\le n\le 2^m \mbox{ and } 1\le k\le d.
\end{equation*}
Then, for 
$s \in [t_{n-1}^m, t_{n}^m )$, we have
\[
d\mathcal{M}^m_s=\mathcal{M}^m_s \left( \omega_{0} (X_{s}^m) ds + \sum_{k=1}^d \omega_{k} (X_{s}^m) \, \Delta B_{{t_{n-1}^{m}}t_{n}^{m}}^{k} ds \right),
\]
Therefore, for $s \in [t_{n-1}^m, t_{n}^m )$, we obtain
\[
\mathcal{M}^m_s=\mathcal{M}^m_{t_{n}^m}e^{- \left( \int_s^{t^m_n}\omega_{0} (X_{u}^m)du  + \sum_{k=1}^d \int_s^{t^m_n}\omega_{k} (X_{u}^m) \Delta B_{{t_{n-1}^{m}}t_{n}^{m}}^{k} du\right)}
\]
Proceeding inductively, we end up with the following identity, valid for $t \in [t_{n-1}^m, t_{n}^m )$ and
$n=0,...,2^m$:
\begin{multline}\label{un}
\mathcal{M}^m_t 
=V(X^m_T)e^{- \left( \int_{t^m_{2^m-1}}^{t^m_{2^m}}\omega_{0} (X_{s}^m)ds  + \sum_{k=1}^d \int_{t^m_{2^m-1}}^{t^m_{2^m}}\omega_{i} (X_{s}^m) \Delta B_{{t_{2^m-1}^{m}}t_{2^m}^{m}}^{k} ds\right)} \times  \\
\cdots \times 
e^{- \left( \int_t^{t^m_n}\omega_{0} (X_{s}^m)ds  + \sum_{k=1}^d \int_t^{t^m_n}\omega_{k} (X_{s}^m) \Delta B_{{t_{n-1}^{m}}t_{n}^{m}}^{k}ds \right)}.
\end{multline}
Owing to the skew-symmetry of $\om_k$ for $k\ge 1$, we can now apply Lemma \ref{contraction} to expression (\ref{un}) in order to get
\begin{align}\label{eq:bnd-M-antisym}
\left\| \mathcal{M}^m_t \right\|  & \le  \exp \left( \frac{T}{2^m} \|  \omega_{0} (X_{t_{2^m-1}^m}^m)  \|
\right)\cdots \exp \left(  (t_{n}^m-t ) \| \omega_{0} (X_{t_{n-1}^m}^m)  \|
\right)\| V(X^m_T)\| \notag \\
&  \le M   \exp \left( CT \right). 
\end{align}
This is our claimed uniform bound on ${M}^m_t$, from which the end of our proof is easily deduced.

\end{proof}

Once the bound (\ref{eq:bnd-first-deriv-infty}) on $\left\| \mathbf{D} X_T^{x} \right\|_{\infty}$ is obtained, one can also retrieve some information on the Hölder norms of $\mathbf{D} X_T^{x}$ improving the general estimate (\ref{eq:bnd-lin-sde}). This is the content of the following proposition:
\begin{proposition}\label{Prop: bnd-first-deriv-holder}
Consider $\frac{1}{2}<\ga<H$ and set $c_0^T=M   \exp( CT)$, which is the constant appearing in relation (\ref{eq:bnd-first-deriv-infty}). Under Assumption~\ref{H}, the Malliavin derivative of the solution $X_T^x$ to equation (\ref{eq:sde}) can be bounded as follows for any $T\in[0,1]$:
\begin{equation}\label{eq:bnd-first-deriv-holder}
\left\| \mathbf{D} X_T^{x} \right\|_{\ga}    \le c_{T,V,d} \, 
\lp 1+|x|+\|B\|_{\ga}^{1/\ga}\rp \, \|B\|_{\ga}^{(1-\ga)/\ga},
\end{equation}
for a strictly positive constant $c_{T,V,d}$.
\end{proposition}
\begin{proof}
We have shown at Theorem \ref{bound MalliavinD} that $\mathbf{D} X_T^{x}$ is governed by equation (\ref{eq:ed-def}), and that $\|\cm\|_{\infty}\le c_0^T$. We will now separate our proof into a local and a global estimate, and notice that the constants appearing in the computations below might change from line to line.

\smallskip

\noindent
\textit{Step1: Local estimate.}
Consider $0\le s<t\le T$ and set $\ep=t-s$. Let $u<v$ be  two generic elements of $[s,t]$. Applying relation (\ref{eq:ineq-young-sharp}) to the expression of $\cm_v-\cm_u$ given by equation (\ref{eq:ed-def}), we obtain
\begin{multline*}
|\cm_v-\cm_u| \le
c_0^T \, c_V |v-u|  \\
+\sum_{i=1}^{d}\lp  
|\cm_u| |\om_i(X_u^x)| \|B\|_{\ga} |v-u|^{\ga}
+\|\cm \om_i(X^x)\|_{s,t,\ga} \|B\|_{\ga} |v-u|^{2\ga}
\rp.
\end{multline*}
We thus obtain, for a constant $c_V$ depending on the vector fields $V$,
\begin{align}\label{eq:local-bnd-M-1}
&\|\cm\|_{s,t,\ga}\le c_0^T c_V |t-s|^{1-\ga} + d \, c_0^T c_V \|B\|_{\ga}
+c_\ga \lc c_0^T c_V \|X^x\|_{\ga}  +c_V \| \cm \|_{s,t,\ga}\rc \|B\|_{\ga} |t-s|^{\ga}  \notag\\
&\le c_0^T c_V \lp T^{1-\ga}+ d  \, \|B\|_{\ga} \rp 
+ d \, c_0^T c_V \|X^x\|_{\ga} \|B\|_{\ga} |t-s|^{\ga}
+ d \, c_V \|B\|_{\ga} |t-s|^{\ga} \| \cm \|_{s,t,\ga}.
\end{align}
Take now $t-s=\ep$ such that $d \, c_V \|B\|_{\ga} \ep^{\ga}=1/2$, namely $\ep=[2d \, c_V \|B\|_{\ga}]^{-1/\ga}$. Recall also that $ \|X^x\|_{\ga}\le c (1+\|B\|_{\ga}^{1/\ga})$ according to \cite{FV-bk}. It is then easily seen that relation (\ref{eq:local-bnd-M-1}) yields
\begin{equation*}
\|\cm\|_{s,t,\ga}\le c_0^T c_{T,V,d} \lp  1+|x|+\|B\|_{\ga}^{1/\ga}\rp:=a_{T,V,d,B},
\end{equation*}
for a strictly positive constant $c_{T,V,d}$.

\smallskip

\noindent
\textit{Step2: Global estimate.}
We consider now $s,t\in[0,T]$ such that $i\ep\le s<(i+1)\ep\le j\ep\le t <(j+1)\ep$, where $\ep$ has been defined at Step 1.
 Set also
$t_i=s$, $t_k=k\ep$ for $i+1\le k\le j$, and $t_{j+1}=t$. 
Then
\begin{eqnarray}\label{eq:global-holder}
\lln \cm_t-\cm_s\rrn
&=& \lln \sum_{k=i}^{j}\cm_{t_{k+1}}-\cm_{t_{k}}\rrn
\le a_{T,V,d,B} \sum_{k=i}^{j} (t_{k+1}-t_{k})^{\ga} \notag \\
&\le& a_{T,V,d,B} (j-i+1)^{1-\ga} (t-s)^\ga, 
\end{eqnarray}
where we have used the fact that $r\mapsto r^\la$ is a concave function.
Note that the indices $i,j$ above satisfy
$(j-i+1)\le 2T/\ep$. Plugging this into the last series of inequalities, 
we end up with our claim (\ref{eq:bnd-first-deriv-holder}).

\end{proof}

\subsection{Log-Sobolev inequality}

In this section, we present some interesting functional inequalities which are usually studied in a Markov setting; namely, the logarithmic-Sobolev inequality and Poincar\'{e} inequality. As we will see, these inequalities become available in our non-Markov case when we have uniform boundedness for the Malliavin derivative of $X_T$ (see Theorem \ref{bound MalliavinD}). 

We start with the following version of logarithmic Sobolev inequality for the law of $X_T$.

\begin{theorem}\label{log-sobolev}
Let $C$ and $M$ be in Theorem \ref{bound MalliavinD}. We have for all $f\in\mathrm{C}^1$ and $T\in[0,1]$,
\begin{align*}
\me f(X_T)^2\ln f(X_T)^2- \left(\me f(X_T)^2\right)\ln \left(\me f(X_T)^2\right)\leq 2M^2e^{2CT}T^{2H}\me|\nabla f(X_T)|^2.
\end{align*}
provided the right hand side in the above is finite.
\end{theorem}
\begin{proof}
The proof is standard by applying Clark-Ocone formula (see e.g. \cite{CHL}). First recall the representation of fractional Brownian motion
$$B_t=\int_0^t K_H(t,s)dW_s.$$
Here $W$ is a $d$-dimensional Wiener process.  Denote by $\mathbf{D}^W$ the Malliavin derivative with respect to the Wiener process $W$. We have
\begin{align}
K_H^*\mathbf{D}=\mathbf{D}^W,
\end{align}
where $K_H^*$ is the isometry from $\mathcal{H}$, the reproducing kernel space of $B$, to $L^2$.
By Clark-Ocone formula we have
\begin{align*}
f(X_T)-\me f(X_T)=\int_0^T\me\big[\mathbf{D}^W_sf(X_T)|\msf_s\big]dW_s=\int_0^T\me\big[K_H^*(\mathbf{D}f(X_T))_s|\msf_s\big]dW_s.
\end{align*}
Hence, if we denote $M_s=\me(f(X_T)|\msf_s), 0\leq s\leq T$, we have
$$dM_s=\me\big[K_H^*(\mathbf{D}f(X_T))_s|\msf_s\big]dW_s.$$
For simplicity we may assume that $f\geq \varepsilon$ for some $\varepsilon>0$, which can be removed afterwards by letting $\varepsilon$ tend to $0$. Applying It\^{o}'s formula to $M_s\ln M_s$, we get 
\begin{eqnarray}\label{M_t ln M_t}
\me f(X_T) \ln(f(X_T)) - \me f(X_T) \ln\lp  \me f(X_T)\rp&=&
\me (M_T\ln M_T) - \me (M_0\ln M_0) \notag\\
&=&\frac{1}{2}\me\int _0^T \frac{1}{M_s}\left|\me\big[K_H^*(\mathbf{D}f(X_T))_s|\msf_s\big]\right|^2ds.
\end{eqnarray}
Replace now $f$ by $f^2$ in the above.  By the Cauchy-Schwarz inequality,
\begin{align*}
\left|\me\big[K_H^*(\mathbf{D}f(X_T)^2)_s|\msf_s\big]\right|^2&=4\left|\me\big[f(X_T)\langle\nabla f(X_T),K_H^*(\mathbf{D}X_T)_s\rangle|\msf_s\big]\right|^2\\
&\leq 4\me\left[f(X_T)^2|\msf_s\right]\me\left[\langle\nabla f(X_T),K_H^*(DX_T)_s\rangle^2|\msf_s\right].
\end{align*}
Substituting the above to (\ref{M_t ln M_t}), together with Theorem \ref{bound MalliavinD}, we obtain the desired result.
\end{proof}

As a corollary of the logarithmic Sobolev inequality obtained above, we have the following Poincar\'{e} inequality (see e.g. \cite[Theorem 8.6.8]{Hsu}). 

\begin{theorem}
Let $C$ and $M$ be in Theorem \ref{bound MalliavinD}. We have for all $f\in \mathrm{C}^1$, 
\begin{align*}
\me f(X_T)^2-\big(\me f(X_T)\big)^2\leq M^2e^{2CT}T^{2H}\me |\nabla f(X_T)|^2,
\end{align*}
provided the right hand side in the above is finite.

\end{theorem}

\begin{remark}
In the above, assume further that the vector fields $V_1,..., V_d$ form an uniform elliptic system, we obtain the following natural expression of logarithmic-Sobolev inequality and Poincar\'{e} inequality when working on a Riemannian manifold
\begin{align*}
\me f(X_T)^2\ln f(X_T)^2- \left(\me f(X_T)^2\right)\ln \left(\me f(X_T)^2\right)\leq 2aM^2e^{2CT}T^{2H}\sum_{i=1}^d\me |V_i f|^2,
\end{align*}
and
\begin{align*}
\me f(X_T)^2-\big(\me f(X_T)\big)^2\leq aM^2e^{2CT}T^{2H}\sum_{i=1}^d\me |V_i f|^2.
\end{align*}
\end{remark}


\subsection{Concentration inequality}

It is classical that the boundedness of the Malliavin derivative in $\mathcal{H}$ implies the Gaussian concentration inequality, more precisely (see \cite[Theorem 9.1.1]{Us} or \cite{Ledoux}):

\begin{lemma}\label{concentration}
Let $F \in \mathbb{D}^{1,2}(\mathcal{H})$ such that almost surely $\| \mathbf{D} F \|_\mathcal{H} \le C$, where $C$ is a non random constant. Then, for every $\theta \ge 0$,
\[
\mathbb{E}(e^{\theta F})  \le e^{ \mathbb{E}(F) \theta + \frac{1}{2} C^2 \theta^2  }
\]
and therefore for every $x \ge 0$,
\begin{equation}\label{eq:concentration-general}
\mathbb{P}\left( F-\mathbb{E}(F) \ge x \right)\le e^{-\frac{x^2}{2C^2}}.
\end{equation}
\end{lemma}

As a corollary of this lemma, we deduce

\begin{proposition}\label{concentration X}
Assume that the Assumption \ref{H} is satisfied, then there exist $C$ and $M$ such that for every $T \ge 0$ and  $\lambda \ge 0$,
\[
\mathbb{P}\left( \sup_{0 \le t \le T} | X^x_t | -\mathbb{E}\left( \sup_{0 \le t \le T} | X^x_t | \right) \ge \lambda \right)\le \exp\left(-\frac{\lambda^2}{2M^2 e^{2CT} T^{2H}}\right).
\]
\end{proposition}

\begin{proof}
Let
$$F=\sup_{0\leq t\leq T}|X_t^x|.$$
By Theorem \ref{bound MalliavinD}, it is not hard to see that (see e.g.  \cite{Nu06})
$$\|\mathbf{D}F\|_\mathcal{H} \leq Me^{CT}T^{H},$$
where $M$ and $C$ are the same as in Theorem \ref{bound MalliavinD}.  Now an easy application of Lemma \ref{concentration} completes our proof. 
\end{proof}

\begin{remark}
Notice that this relation can only be obtained for $H>\frac{1}{2}$. Indeed, denote by $\cac^0$ the space of continuous functions endowed with the supremum norm. Then  the norm involved in \cite[Theorem 9.1.1]{Us} is $\|\mathbf{D}X_t\|_{L^\infty(\ch)}$. This norm is dominated by $\|\mathbf{D}X_t\|_{L^\infty(\cac^{0})}$ only if $H>\frac{1}{2}$.
\end{remark}

\subsection{Gaussian upper bound}\label{sec:gauss-bnd}
One natural way to estimate the density of $ X_t$ is to apply the results in \cite[Chapter 2]{Nu06}. More precisely, we first have the following integration by parts formula for non-degenerate random vectors.
\begin{proposition}
Let $F=(F^1,...,F^d)$ be a non-degenerate random vector and $\gamma_F$ the Malliavin matrix of $F$. Let $G\in\D^\infty$ and $\varphi$ be a function in the space $C_p^\infty(\mr^d)$. Then for any multi-index $\alpha\in\{1,2,...,d\}^k, k\geq 1$, there exists an element $H_\alpha\in\D^\infty$ such that
$$\me[\partial_\alpha \varphi(F)G]=\me[\varphi(F)H_\alpha].$$
Moreover, the elements $H_\alpha$ are recursively given by
\begin{align*}
&H_{(i)}=\sum_{j=1}^{d}\delta\left(G(\gamma^{-1})^{ij}\mathbf{D}F^j\right)\\
&H_\alpha=H_{\alpha_k}(H_{(\alpha_1,..., \alpha_{k-1})}),
\end{align*}
and for $1\leq p<q<\infty$ we have
$$\|H_\alpha\|_{L^p}\leq C_{p,q}\|\gamma_F^{-1}\mathbf{D}F\|^k_{k, 2^{k-1}r}\|G\|_{k,q},$$
where $\frac{1}{p}=\frac{1}{q}+\frac{1}{r}$.
\end{proposition}
\begin{remark}
By the estimates for $H_\alpha$ above, one can conclude that there exist constants $\beta, \gamma>1$ and integers $m, n$ such that
\begin{align}\label{est H}
\|H_\alpha\|_{L^p}\leq C_{p,q}\|\det\gamma^{-1}_F\|^m_{L^\beta}\|\mathbf{D}F\|_{k,\gamma}^n\|G\|_{k,q}.
\end{align}
\end{remark}
\begin{remark}
In what follows, we use $H_\alpha(F,G)$ to emphasize its dependence on $F$ and $G$.
\end{remark}

As a consequence of the above proposition, one has
\begin{proposition}\label{rep density}
Let $F=(F^1,..., F^d)$ be a non-degenerate random vector. Then the density $p_F(x)$ of $F$ belongs to the Schwartz space, and
\begin{align*}
p_F(x)=\me[ \mathrm{I}_{\{F>x\}}H_{(1,2,...,d)}(F,1)].
\end{align*}
\end{proposition}

Now we state and prove a global Gaussian upper bound for the density function of $X^x_t$.

\begin{theorem}\label{thm:main-gauss-bnd}
Denote by $p_X(t,y)$ the density function of $X_t^x$. There exist positive constants $c^{(1)}_t, c^{(2)}_{t,x}, c_3, \mathrm{and}\  c_4$ such that for all $t\in[0,\infty)$,
\begin{align}
p_X(t,y)\leq c^{(1)}_t \exp\left(-\frac{\left(|y|-c^{(2)}_{t,x}\right)^2}{c_3 e^{c_4t} t^{2H}}\right).
\end{align}
Here $c^{(1)}_t$ is of the form $$c^{(1)}_t= O\left(\frac{1}{t^\alpha}\right)\quad\quad\quad\mathrm{as}\ t\downarrow 0,$$ for some positive number $\alpha$; and $c^{(2)}_{t,x}$ converges to a constant as $t\downarrow 0$.
\end{theorem}

\begin{proof}
Fix $\beta\in (\frac{1}{2}, 1)$.  By (\ref{est H}) and Proposition \ref{rep density}, we have
\begin{align}\label{rep-density-final}
p_X(t,y)
\leq C(\mp\{X_t>y\})^\frac{1}{q} \|\det\gamma^{-1}_{X_t}\|^m_{L^\alpha}\|\mathbf{D}X_t\|_{k,\gamma}^n,
\end{align}
for some constants $\alpha, \gamma>1$ and integers $m, n$. Without lose of generality, we may assume $y_i\geq0$ for $1\le i\le d$. Since otherwise, for example $y_j<0$, we can consider the alternative expression for the density
$$p_X(t,y)=\me\left[\Pi_{i\not = j}I_{\{y^i<X^i_t\}}I_{\{y^j>X^j_t\}}H_{(1,2,...,d)}(X_t,1)\right],$$
and deduce similar estimate. In the following, we provide estimate for $\|\mathbf{D}X_t\|_{k,\gamma}^n$, $\|\det\gamma_{X_t}^{-1}\|_{L^\alpha}^m$ and $\mp\{X_t>y\}$ respectively.

By Proposition \ref{prop:moments-sdes} there exist constants $C >0$, depending on $V, x$ and $k$, such that
\begin{align}\label{X-est}
&\sup_{0\leq t\leq T}|X_t^i|\leq |x^i|+CT\|B\|_{0,T,\beta}^{1/\beta},\\ \label{gamma-est}
&\sup_{0\leq t\leq T}\|\gamma_{X_t}^{-1}\|\leq \frac{C}{T^{2Hd}}\left[1+e^{CT\|B\|_{0,T,\beta}^{1/\beta}}\right],\\ \label{DX-est}
&\sup_{0\leq t, r_i\leq T}|\mathbf{D}_{r_k}^{j_k}...\mathbf{D}_{r_1}^{j_1}X_t^i|\leq Ce^{CT\|B\|_{0,T,\beta}^{1/\beta}}.
\end{align}
On the other hand we have, for some constant $M_\beta$ (cf. \cite{HN}):
\begin{align}\label{tail B}
\mp\left\{\|B\|_{0,T,\beta}>r\right\}\leq M_\beta e^{-\frac{r^2}{2T^{2(H-\beta)}}}.
\end{align}
Hence
\begin{align*}
&\|\det\gamma_{X_t}^{-1}\|_{L^\alpha}^m<\infty\quad\mathrm{and}\quad\|\mathbf{D}X_t\|^n_{k,\gamma}<\infty.
\end{align*}
Moreover,  by (\ref{gamma-est}), (\ref{DX-est}) and the tail estimate (\ref{tail B}), there exists constant $C>0$ independent of $t$ such that
\begin{align*}
&\|\det\gamma_{X_t}^{-1}\|_{L^\alpha}^m\leq \frac{C}{t^{2Hmd}}\left(1+t^\frac{m}{\alpha}F(t,1/\beta-1,1/\beta,C)^{\frac{m}{\alpha}}\right),
\end{align*}
and
\begin{align*}
&\|\mathbf{D}X_t\|^n_{k,\gamma}\leq C\left(1+t^{nH(k+1))+\frac{n}{\gamma}}F(t,1/\beta-1,1/\beta,C)^{\frac{n}{\gamma}}\right).  
\end{align*}
In the above
\begin{align*}
F(t,a,b,M)=\int_0^\infty r^{a}e^{-\frac{r^2}{2t^{2(H-\beta)}}}e^{Mtr^{b}}dr,
\end{align*}
for any $a, M>0$ and $2>b>0$. Now set $$c^{(1)}_t=\|\det\gamma_{X_t}^{-1}\|_{L^\alpha}^m\|\mathbf{D}X_t\|^n_{k,\gamma}.$$ Elementary computation shows, for some constant $\alpha>0$
\begin{align}\label{c^1_t}c^{(1)}_t\sim O\left(\frac{1}{t^\alpha}\right) \quad\quad\mathrm{as}\ t\downarrow 0.\end{align}

Finally we estimate $\mp\{X_t>y\}$. Define 
$$\xi_t=\max_{1\leq i\leq d}\left\{|X^1_t|, |X^2_t|,...,|X^d_t|\right\}.$$
First note that $$|\mathbf{D}\xi_t|\leq Me^{Ct},$$
since the Malliavin derivative $|\mathbf{D}X^i_t|\leq Me^{Ct}$ for each $i=1,2,...d$. Hence $\xi_t$ has the same concentration property as specified in Proposition \ref{concentration X}.  

By (\ref{X-est}) and the concentration property for $\zeta_t$, we conclude 
\begin{align}\label{est tail}
\mp\left\{X_t>y\right\}\leq \mp\left\{ \xi_t>\frac{|y|}{\sqrt{d}}\right\}
\leq\exp\left(-\frac{\left(\frac{|y|}{\sqrt{d}}-\me \xi_t\right)^2}{2M e^{Ct} t^{2H}}\right)
\leq\exp\left(-\frac{\left({|y|}-c^{(2)}_{t,x}\right)^2}{2\sqrt{d}M e^{Ct} t^{2H}}\right).
\end{align}
Now (\ref{rep-density-final}), (\ref{c^1_t}) and (\ref{est tail}) give us the desired upper bound for the density $p_X(t,y)$.
\end{proof}

\begin{remark} 
There is also an upper bound for the constant $c^{(1)}_t$ in the above theorem as $t\uparrow \infty$. Indeed, by some elementary computation one can show that 
$$c^{(1)}_t\leq t^\alpha e^{t^\beta}\quad\mathrm{as}\ t\uparrow \infty,$$
for some $\alpha, \beta>0$.
\end{remark}

\section{Extension to the irregular case}
\label{sec:irregular}

From now on, our purpose is to extend the previous to the case of a fBm with Hurst index $\frac{1}{3}<H<\frac{1}{2}$. This requires the introduction of some rough paths tools, which is the aim of the current section. We shall use in fact the language of algebraic integration theory, which is a variant of the rough paths theory introduced in~\cite{Gu} (we also refer to \cite{GT} for a detailed introduction of the topic).

\subsection{Increments}\label{incr}

The extended pathwise integration we will deal with is based on the
notion of \emph{increments}, together with an elementary operator $\der$
acting on them. The algebraic structure they generate is described
in \cite{Gu,GT}, but here we present  directly the definitions of
interest for us, for sake of conciseness. First of all,  for an
arbitrary real number $T>0$, a vector space $V$ and an integer $k\ge
1$ we denote by $\cac_k(V)$ the set of functions $g : [0,T]^{k} \to
V$ such that $g_{t_1 \cdots t_{k}} = 0$ whenever $t_i = t_{i+1}$ for
some $i\le k-1$. Such a function will be called a
$(k-1)$-increment, and we   set $\cac_*(V)=\cup_{k\ge
1}\cac_k(V)$. We can now define the announced elementary operator
$\der$ on $\cac_k(V)$:
\begin{equation}
 \label{eq:coboundary}
\delta : \cac_k(V) \to \cac_{k+1}(V), \qquad
(\delta g)_{t_1 \cdots t_{k+1}} = \sum_{i=1}^{k+1} (-1)^{k-i}
g_{t_1  \cdots \hat t_i \cdots t_{k+1}} ,
\end{equation}
where $\hat t_i$ means that this particular argument is omitted.
A fundamental property of $\der$, which is easily verified,
is that
$\delta \delta = 0$, where $\delta \delta$ is considered as an operator
from $\cac_k(V)$ to $\cac_{k+2}(V)$.
We    denote $\cz\cac_k(V) = \cac_k(V) \cap \text{Ker}\delta$
and $\cb \cac_k(V) =
\cac_k(V) \cap \text{Im}\delta$.

\vspace{0.3cm}

Some simple examples of actions of $\der$,
which will be the ones we will really use throughout the paper,
are obtained by letting
$g\in\cac_1$ and $h\in\cac_2$. Then, for any $s,u,t\in\ott$, we have
\begin{equation}
\label{eq:simple_application}
 (\der g)_{st} = g_t - g_s,
\quad\mbox{ and }\quad
(\der h)_{sut} = h_{st}-h_{su}-h_{ut}.
\end{equation}
Furthermore, it is easily checked that
$\cz \cac_{k}(V) = \cb \cac_{k}(V)$ for any $k\ge 1$.
In particular, the following basic property holds:
\begin{lemma}\label{exd}
Let $k\ge 1$ and $h\in \cz\cac_{k+1}(V)$. Then there exists a (non unique)
$f\in\cac_{k}(V)$ such that $h=\der f$.
\end{lemma}

\begin{proof} 
This elementary proof is included in
\cite{Gu}, and will be omitted here. However, let us  mention that
$f_{t_1\ldots t_{k}}=(-1)^{k+1}h_{0t_1\ldots t_{k}}$ is a possible
choice. 

\end{proof}

Observe that Lemma \ref{exd} implies that all the elements
$h \in\cac_2(V)$ such that $\der h= 0$ can be written as $h = \der f$
for some (non unique) $f \in \cac_1(V)$. Thus we get a heuristic
interpretation of $\der |_{\cac_2(V)}$:  it measures how much a
given 1-increment  is far from being an  exact increment of a
function, i.e., a finite difference.

\vspace{0.3cm}

Notice that our future discussions will mainly rely on
$k$-increments with $k \le 2$, for which we will make  some
analytical assumptions. Namely,
we measure the size of these increments by H\"older norms
defined in the following way: for $f \in \cac_2(V)$ let
\begin{equation}\label{eq:def-norm-C2}
\norm{f}_{\mu} =
\sup_{s,t\in\ott}\frac{|f_{st}|}{|t-s|^\mu},
\quad\mbox{and}\quad
\cac_2^\mu(V)=\lcl f \in \cac_2(V);\, \norm{f}_{\mu}<\infty  \rcl.
\end{equation}
Obviously, the usual H\"older spaces $\cac_1^\mu(V)$ will be determined
       in the following way: for a continuous function $g\in\cac_1(V)$, we simply set
\begin{equation}\label{def:hnorm-c1}
\|g\|_{\mu}=\|\der g\|_{\mu},
\end{equation}
and we will say that $g\in\cac_1^\mu(V)$ iff $\|g\|_{\mu}$ is finite.
Notice that $\|\cdot\|_{\mu}$ is only a semi-norm on $\cac_1(V)$,
but we will generally work on spaces of the type
\begin{equation}\label{def:hold-init}
\cac_{1,a}^\mu(V)=
\lcl g:\ott\to V;\, g_0=a,\, \|g\|_{\mu}<\infty \rcl,
\end{equation}
for a given $a\in V$, on which $\|g\|_{\mu}$ defines a distance in
the usual way.
For $h \in \cac_3(V)$ set in the same way
\begin{eqnarray}
 \label{eq:normOCC2}
 \norm{h}_{\gamma,\rho} &=& \sup_{s,u,t\in\ott}
\frac{|h_{sut}|}{|u-s|^\gamma |t-u|^\rho}\\
\|h\|_\mu &= &
\inf\left \{\sum_i \|h_i\|_{\rho_i,\mu-\rho_i} ;\, h =
\sum_i h_i,\, 0 < \rho_i < \mu \right\} ,\nonumber
\end{eqnarray}
where the last infimum is taken over all sequences $\{h_i \in \cac_3(V) \}$
such that $h
= \sum_i h_i$ and for all choices of the numbers $\rho_i \in (0,\mu)$.
Then  $\|\cdot\|_\mu$ is easily seen to be a norm on $\cac_3(V)$, and we set
$$
\cac_3^\mu(V):=\lcl h\in\cac_3(V);\, \|h\|_\mu<\infty \rcl.
$$
Eventually,
let $\cac_3^{1+}(V) = \cup_{\mu > 1} \cac_3^\mu(V)$,
and notice  that the same kind of norms can be considered on the
spaces $\cz \cac_3(V)$, leading to the definition of some spaces
$\cz \cac_3^\mu(V)$ and $\cz \cac_3^{1+}(V)$.

\vspace{0.3cm}

With these notations in mind
the following proposition is a basic result, which  belongs to  the core of
our approach to pathwise integration. Its proof may be found
in a simple form in \cite{GT}.
\begin{proposition}[The $\Lambda$-map]
\label{prop:Lambda}
There exists a unique linear map $\Lambda: \cz \cac^{1+}_3(V)
\to \cac_2^{1+}(V)$ such that
$$
\delta \Lambda  = \id_{\cz \cac_3^{1+}(V)}
\quad \mbox{ and } \quad \quad
\Lambda  \delta= \id_{\cac_2^{1+}(V)}.
$$
In other words, for any $h\in\cac^{1+}_3(V)$ such that $\der h=0$
there exists a unique $g=\laa(h)\in\cac_2^{1+}(V)$ such that $\der g=h$.
Furthermore, for any $\mu > 1$,
the map $\laa$ is continuous from $\cz \cac^{\mu}_3(V)$
to $\cac_2^{\mu}(V)$ and we have
\begin{equation}\label{ineqla}
\|\Lambda h\|_{\mu} \le \frac{1}{2^\mu-2} \|h\|_{\mu} ,\qquad h \in
\cz \cac^{\mu}_3(V).
\end{equation}
\end{proposition}

\smallskip

Let us mention at this point a first link between the structures
we have introduced so far and the problem of integration of irregular
functions.

\begin{corollary}
\label{cor:integration}
For any 1-increment $g\in\cac_2 (V)$ such that $\der g\in\cac_3^{1+}$,
set
$
\delta f = (\id-\Lambda \delta) g
$.
Then
$$
(\delta f)_{st} = \lim_{|\Pi_{st}| \to 0} \sum_{i=0}^{n-1}
g_{t_i\,t_{i+1}},
$$
where the limit is over any partition $\Pi_{st} = \{t_0=s,\dots,
t_n=t\}$ of $[s,t]$, whose mesh tends to zero. Thus, the 1-increment
$\delta f$ is the indefinite integral of the 1-increment $g$.
\end{corollary}

\subsection{Computations in $\cac_*$}\label{cpss}

Let us specialize now to the case $V=\R$, and just write $\cac_{k}^{\ga}$ for $\cac_{k}^{\ga}(\R)$. Then $(\cac_*,\delta)$ can be endowed with the following product:
for  $g\in\cac_n$ and $h\in\cac_m$ let  $gh$
be the element of $\cac_{n+m-1}$ defined by
\begin{equation}\label{cvpdt}
(gh)_{t_1,\dots,t_{m+n+1}}=
g_{t_1,\dots,t_{n}} h_{t_{n},\dots,t_{m+n-1}},
\quad
t_1,\dots,t_{m+n-1}\in\ott.
\end{equation}
In this context, we have the following useful properties.

\begin{proposition}\label{difrul}
The following differentiation rules hold true:
\begin{enumerate}
\item
Let $g\in\cac_1$ and $h\in\cac_1$. Then
$gh\in\cac_1$ and $\der (gh) = \der g\,  h + g\, \der h$.

\item
Let $g\in\cac_1$ and $h\in\cac_2$. Then
$gh\in\cac_2$ and $\der (gh) =- \der g\, h + g \,\der h$.

\item
Let $g\in\cac_2$ and $h\in\cac_1$. Then
$gh\in\cac_2$ and $\der (gh) = \der g\, h  + g \,\der h$.

\end{enumerate}
\end{proposition}

\smallskip

The iterated integrals of smooth functions on $\ott$ are obviously
particular cases of elements of $\cac$, which will be of interest
for us. Let us recall  some basic  rules for these objects: consider
$f\in\cac_1^\infty$ and $g\in\cac_1^\infty$, where $\cac_1^\infty $
denotes the set of smooth functions on $\ott$. Then the integral
$\int f \, dg$, which will be denoted  indistinctly by $\int f\, dg$
or $\cj(f\, dg)$, can be considered as an element of
$\cac_2^\infty$. Namely, let $\cs_{2,T}$ denote the simplex $\{(s,t)\in [0,T]^2: 0\leq s< t\leq T\}$, for $(s,t)\in\cs_{2,T}$ we set
$$
\cj_{st}(f\,  dg)
=
\left(\int  f dg \right)_{st} = \int_s^t   f_u dg_u.
$$
The multiple integrals can also be defined in the following way:
given a smooth element $h \in \cac_2^\infty$ and $(s,t)\in\cs_{2,T}$, we set
$$
\cj_{st}(h\, dg )\equiv
\left(\int h dg  \right)_{st} = \int_s^t  h_{su} dg_u .
$$
In particular,  for
$f^1\in\cac_1^\infty$, $f^2\in\cac_1^\infty$
and $f^3\in\cac_1^\infty$ the double integral
$\cj_{st}( f^3\, df^2 df^1)$ is defined  as
$$
\cj_{st}( f^3\, df^2df^1)
=\lp \int f^3\, df^2 df^1  \rp_{st}
= \int_s^t \cj_{su}\lp f^3\, df^2  \rp \, df_u^1.
$$
Now suppose that the $n$th order iterated integral of
$f^{n+1}df^n\cdots df^2$, which is  denoted by
$\cj(f^{n+1}df^n\cdots df^2)$, has been defined for
$f^j\in\cac_1^\infty$.
Then, if $f^1\in\cac_1^\infty$, we set
\begin{equation}\label{multintg}
\cj_{st}(f^{n+1}df^n \cdots df^2 df^1)
=
\int_s^t  \cj_{su}\lp f^{n+1}df^n\cdots df^2\rp  \, df_u^{1},
\end{equation}
which recursively defines the iterated integrals of smooth
functions. Observe that an   $n$th order integral $\cj(df^n\cdots df^2
df^1)$ can be defined along the same lines, starting with
\begin{equation*}
\cj(df)=\delta f, \qquad
\cj_{st}(df^2\,df^1)=\int_s^t \cj_{su}(df^2)\,df^1_u=\int_{s}^t \big( \delta
f^2\big)_{su}\,df^1_u,
\end{equation*}
and so on.

\smallskip

The following relations between multiple integrals and the operator $\der$ will also be useful. The reader is sent to \cite{GT} for its elementary proof.
\begin{proposition}\label{dissec}
Let $f\in\cac_1^\infty$ and $g\in\cac_1^\infty$.
Then it holds that
$$
\der g = \cj( dg), \qquad \der\lp \cj(f dg)\rp = 0, \qquad \der\lp
\cj (df dg)\rp = (\der f) (\der g) = \cj(df) \cj(dg),
$$
and
$$
\der \lp \cj( df^n \cdots df^1)\rp  =
\sum_{i=1}^{n-1}
\cj\lp df^n \cdots df^{i+1}\rp \cj\lp df^{i}\cdots df^1\rp.
$$
\end{proposition}

\subsection{Weakly controlled processes}\label{sec:wc-ps}

The rough path theory allows to define and solve differential equations driven by a generic Hölder continuous path $B$ provided enough iterated integrals of this function can be defined. We shall briefly recall how this is done, in the simplest nontrivial case of a Hölder continuity exponent $\frac{1}{3}<\ga<\frac{1}{2}$. Observe that we keep here the notation $B$ for the underlying path as in the fBm case for notational sake, while the theory can be applied to much more general situations.

\smallskip

The basic assumption one has to add in order to define our objects when $\ga>\frac{1}{3}$ is the existence of an (abstract) double iterated integral of $B$ with respect to itself, which can be defined as follows:
\begin{hypothesis}\label{hyp:x}
The path $B$ is $\R^d$-valued $\ga$-H\"older with $\ga>\frac{1}{3}$ and
admits a L\'evy area, that is a process
$\xd\in\cac_2^{2\ga}(\R^{d,d})$ satisfying
$$\,
\der\xd=\bb^{\1}\otimes \bb^{\1},
\quad\mbox{i.\!\! e.}\quad
\bb^{\2,ij}_{sut} 
=
[\bb^{\1,i}]_{su} [\bb^{\1,j}]_{ut},
$$
for $s,u,t\in\cs_{3,T}$ and $i,j\in\{1,\ldots,d  \}$. We also assume that $\bb^{\2}$ can be obtained in the following way: consider the sequence of linear dyadic approximation $B^m$ of $B$ defined like in  Proposition \ref{prop:lin-interpol}. For $0\le s<t\le 1$ and $i_1,i_2\in\{1,\ldots,d\}$, set $\bb^{\2,m,i_1i_2}_{st}=\int_{s<u_1<u_2<t} dB_{u_1}^{m,i_1} \, dB_{u_2}^{m,i_2}$, which is defined as a Riemann-Stieljes integral. Then we suppose that $\bb^{\2,m}$ converges to $\bb^{\2}$ in the norm of $\cac_2^{2\ga}$.
\end{hypothesis}
It should be noticed at this point that fBm satisfies the above assumption, as shown in~\cite{CQ,FV-bk,NTU}:

\begin{proposition}
Let $B$ be a $d$-dimensional fBm with $H>\frac{1}{3}$ as defined in Section \ref{sec:stoch-calc-fbm}. For the linear dyadic approximation $B^m$ of $B$ defined at Proposition \ref{prop:lin-interpol}, the increment $\bb^{\2,m}$ almost surely converges to an  element $\bb^{\2}$ satisfying Hypothesis \ref{hyp:x}. The convergence holds in any $\cac_2^{2\ga}$ norm for $\ga<H$.
\end{proposition}

\smallskip

The first difference between the Young case and the situation of a Hölder continuity exponent $\frac{1}{3}<\ga\le \frac{1}{2}$ is that a restriction has to be imposed on the class of allowed integrands with respect to $B$. This class is called the class of \emph{weakly controlled processes}, and is defined as follows:
\begin{definition}\label{def:weak-ctrled-ps}
Let $z$ be a process in $\cac_1^\ga(\R^n)$ with $\frac{1}{3}<\ga\le \frac{1}{2}$
(that is, $N:=\lfloor 1/\ga\rfloor=2$). We say that $z$ is a weakly controlled path based
on $B$ and starting from $a$ if $z_0=a$, which is a given initial condition in $\R^n$,
and $\der z\in\cac_2^\gamma(\R^n)$ can be decomposed into
\begin{equation}\label{weak:dcp}
\der z^{i}=\zeta^{i i_1} \bb^{\1,i_1}+ r^{i},
\quad\mbox{i.\!\! e.}\quad
(\der z^{i})_{st}=\zeta_s^{i i_1} \bb^{\1,i_1}_{st} + r_{st}^{i},
\end{equation}
for all $(s,t)\in\cs_{2,T}$. In the previous formula, we assume
$\zeta\in\cac_1^\ga(\R^{n,d})$, and $r$ is a regular part
such that $r\in\cac_2^{2\ga}(\R^n)$. The space of weakly controlled
paths starting from $a$ will be denoted by $\cq_{\ga,a}(\R^n)$, and a process
$z\in\cq_{\ga,a}(\R^n)$ can be considered in fact as a couple
$(z,\zeta)$. The natural semi-norm on $\cq_{\gamma,a}(\R^n)$ is given
by
$$
\cn[z;\cq_{\ga,a}(\R^n)]=
\cn[z;\cac_1^{\ga}(\R^n)]
+ \cn[\zeta;\cac_1^{\infty}(\R^{n,d})]
+ \cn[\zeta;\cac_1^{\ga}(\R^{n,d})]
+\cn[r;\cac_2^{2\ga}(\R^n)],
$$
with $\cn[g;\cac_1^{\ka}]$ defined by
(\ref{def:hnorm-c1}) and
$\cn[\zeta;\cac_1^{\infty}(V)]=\sup_{0\le s\le T}|\zeta_s|_V$.
\end{definition}

Two basic steps in order to define and solve differential equations with respect to $B$ are then:
\begin{enumerate}
\item
Study the decomposition of $f(z)$ as weakly controlled process, when $f$ is a smooth function and $z$ a weakly controlled process.
\item
Define rigorously the integral $\int z_u dB_u=\cj(z dB)$
for a weakly controlled path $z$ and compute its decomposition
(\ref{weak:dcp}).
\end{enumerate}
We shall now detail a little this program.

\smallskip

Let us see then how to decompose $f(z)$ as a controlled process when $f$ is a smooth enough function, a step for which we first introduce a convention which will hold true until the end of the paper: for any smooth function $f:\R^n\to\R$, $k\ge 1$, $(i_1,\ldots,i_k)\in\{1,\ldots, d\}^k$ and $\xi\in\R^n$, we set
\begin{equation}\label{eq:convention-partial-deriv}
\partial_{i_1\ldots i_k}^{k} f(\xi)=\frac{\partial^k f}{\partial x_{i_1}\cdots \partial x_{i_k}}(\xi).
\end{equation}
With this notation in hand, our decomposition result is the following:
\begin{proposition}\label{cp:weak-phi}
Let $f:\R^n\to\R$ be a $C_b^2$ function such that $f(a)=\ha$, $z$ a controlled process as in Definition \ref{def:weak-ctrled-ps} and set $\hz=f(z)$. Then $\hz\in\cq_{\ga,\ha}(\R)$,  and it can be decomposed into
$\der \hz=\hze^{i_1} \bb^{\1,i_1} +\hr$,
with
$$
\hze^{i_1}= \partial_{i}f(z) \, \ze^{i i_1} 
\quad\mbox{ and }\quad
\hat{r}= \lc \der f(z)- \partial_{i}f(z) \, \der z^{i} \rc + \partial_{i}f(z) \, r^{i}.
$$
Furthermore,
\begin{equation}\label{bnd:phi}
\cn[\hat{z};\cq_{\gamma,\hat{a}}(\mr)]\le c_{f,T}\lp 1+\cn^2[z;\cq_{\ga,a}(\R^n)]
\rp.
\end{equation}
\end{proposition}

\smallskip

Let us now turn  to the integration of weakly controlled paths,
which is summarized in the following theorem.
\begin{theorem}\label{intg:mdx}
For a given $\frac{1}{3}<\ga\le \frac{1}{2}$, let $B$ be a process satisfying
Hypothesis \ref{hyp:x}. Furthermore,  let $m\in\cq_{\ga,b}(\R^{d})$
with $m_0=b\in\R^d$ and decomposition
\begin{equation}\label{dcp:m}
\der m^{i}=\mu^{i i_1} \, \bb^{\1,i_1}+ r^{i},
\quad\mbox{ where }\quad
\mu\in\cac_1^\ga(\R^{d, d}), \, r\in\cac_2^{2\ga}(\R^{d}).
\end{equation}
Define $z$ by $z_0=a\in\R$ and
\begin{equation}\label{dcp:mdx}
\der z=
m^{i} \, \bb^{\1,i} + \mu^{i i_1} \, \bb^{\2,i_1 i}
-\laa\lp  r^{i} \, \bb^{\1,i} + \der\mu^{i i_1} \, \bb^{\2,i_1 i} \rp.
\end{equation}
Finally, set
$$
\cj_{st}(m\, dB)=\ist \left\langle m_u, \, dB_u\right\rangle_{\R^d} \triangleq \der z_{st}.
$$
Then:

\smallskip

\noindent \emph{(1)} $z$ is well-defined as an element of
$\cq_{\ga,a}(\R)$, and coincides with the Riemann-Stieltjes integral
of $m$ with respect to $B$ whenever these two functions are smooth.

\smallskip

\noindent \emph{(2)} The semi-norm of $z$ in $\cq_{\ga,a}(\R)$ can
be estimated as
\begin{equation}\label{bnd:norm-imdx}
\cn[z;\cq_{\ga,a}(\R)]\le c_{B} \lp 1 +
\cn[m;\cq_{\ga,b}(\R^{d})]\rp,
\end{equation}
for a positive constant $c_{B}$ which can be bounded as follows:
$$
c_B\le c
\lp
\cn[\bb^{\1};\, \cac_2^{\ga}(\R^d)]+ \cn[\bb^{\2};\, \cac_2^{2\ga}(\R^{d^2})]
\rp,
\quad\mbox{ for a universal constant }c.
$$

\smallskip

\noindent
\emph{(3)}
It holds
\begin{equation}\label{rsums:imdx}
\cj_{st}(m\, dx) =\lim_{|\Pi_{st}|\to 0}\sum_{q=0}^{n-1} \lc
m_{t_{q}}^{i} \, \bb^{\1}_{t_{q}, t_{q+1}}(i) + \mu_{t_{q}}^{i i_1} \,
\bb^{\2,i_1 i}_{t_{q}, t_{q+1}} \rc,
\end{equation} for any $0\le s<t\le T$,
where the limit is taken over all partitions
$\Pi_{st} = \{s=t_0,\dots,t_n=t\}$
of $[s,t]$, as the mesh of the partition goes to zero.
\end{theorem}

\subsection{Rough differential equations}\label{sec:rdes}

Recall that we are concerned with equations of the form (\ref{eq:sde-intro}).
In our algebraic setting, we will rephrase this as follows: we will
say that $X^x$ is a solution to (\ref{eq:sde-intro}) if $X_0^x=x$,
$X^x\in\cq_{\ga,x}(\R^d)$ and for any $0\le s\le t\le 1$ we have
\begin{equation}\label{eds:alg-form}
\der X^x_{st}=\int_s^t V_0 (X^x_u)du+\cj_{st}\lp V_i(X^x)\, dB^{i}\rp,
\end{equation}
where the integral $\cj(V_i(X^x)\, dB^{i})$ has to be understood in the sense
of Proposition \ref{intg:mdx}. The following existence and uniqueness result is then classical in rough paths theory. 

\begin{theorem}\label{thm:ex-uniq}
Let $B$ be a process satisfying Hypothesis \ref{hyp:x} and $V_0,\ldots,V_n$ a collection of vector fields which fulfill Hypothesis \ref{H}. Then 

\smallskip

\noindent
\emph{(i)}
Equation (\ref{eds:alg-form}) admits a unique solution $X^x$ in
$\cq_{\ga,x}(\R^d)$.

\smallskip

\noindent
\emph{(ii)}
Consider the linear approximation $B^m$ of $B$ introduced in Hypothesis \ref{hyp:x}, and set $\tilde X^m$ for the solution of the (ordinary) differential equation  (\ref{eds:alg-form}) driven by the piecewise smooth function $B^m$. Then $\tilde X^m$ converges to $X^x$ in $\cac_1^\ga$ norm.

\smallskip

\noindent
\emph{(iii)}
Consider the sequence of dyadic partitions of Proposition \ref{prop:lin-interpol}. Define a process $X^{n}$ on the points $t_k^n$ by $X^{n}_0=x$ and 
\begin{equation}\label{eq:davie-scheme}
\der X^{n}_{t_k t_{k+1}}= \frac{V_0(X^{n}_{t_k})}{2^n} + V_{i}(X^{n}_{t_k}) \bb^{\1,i}_{t_k t_{k+1}}
+V_{i_1}V_{i_2}(X^{n}_{t_k}) \bb^{\2,i_1i_2}_{t_k t_{k+1}}.
\end{equation}
Complete the definition of $X^n$ on $\ou$ by linear interpolation. Then as $n\to\infty$, the process $X^n$ converges to $X^x$ in $\cac_1^{\ga}$ norm.
\end{theorem}

\begin{proof}
We refer to \cite{LQ} for the proof of the existence and uniqueness part, as well as to \cite{Gu} for the same result in the algebraic integration setting. Part (ii) of our proposition stems from the continuity of the Itô map, which is also stated and proved in both \cite{LQ,Gu}. The approximation statement (iii) has first been stated by Davie \cite{Da} and then been generalized in \cite{FV-bk}.

\smallskip

In the sequel we will simply try to relate the decomposition of the solution to Equation~(\ref{eds:alg-form}) as a controlled process and the numerical scheme given by (\ref{eq:davie-scheme}), a relation which turns out to be useful in the sequel. For this, we shall denote by $r$ any  increment in $\cac_2^{2\ga}$ and by $r^{\sharp}$ any  increment in $\cac_2^{1^{+}}$ in the computations below, independently of their values. Observe then that, according to the right hand side of (\ref{eds:alg-form}), the decomposition of $X^x$ as a controlled process is given by
\begin{equation*}
\der X^{x,j}=\zeta^{j j_1} \bb^{\1,j_1} + r,
\quad\mbox{with}\quad 
\zeta^{j j_1}=V_{j_1}^{j}(X^x).
\end{equation*}
Hence, owing to Proposition \ref{cp:weak-phi}, one has
\begin{equation}\label{eq:dcp-V-i-X}
\der V_i(X^x)= \hze^{i j_1} \bb^{\1,j_1} + r^{i},
\quad\mbox{with}\quad 
\hze^{i j_1}= \partial_j V_i (X^x)  \ze^{j j_1}=V_{j_1} V_{i}  (X^x).
\end{equation}

Now, if one desires an expansion of $\der X^{x}$ up to increments of regularity $1^+$, consider again the right hand side of (\ref{eds:alg-form}), and compute it by a direct application of Theorem~\ref{intg:mdx}. This yields 
\begin{eqnarray*}
\der X^{x}_{st}&=&V_0(X^x_s) (t-s) + V_i(X_s^x) \bb^{\1,i}_{st} 
+ \hze^{i j_1}  \bb^{\2,j_1 i}_{st} + r^{\sharp} \\
&=&V_0(X^x_s) (t-s) + V_i(X_s^x) \bb^{\1,i}_{st} + V_{j_1} V_{i}  (X^x_s) \bb^{\2,j_1 i}_{st} + r^{\sharp}.
\end{eqnarray*}
Thanks to identity (\ref{rsums:imdx}), it is now easily seen that (\ref{eq:davie-scheme}) is a natural candidate for our numerical scheme.

\end{proof}

\smallskip

We show now how to get efficient bounds on the solution to equation (\ref{eds:alg-form}) out of its numerical scheme. This step is understood as a warmup for the same kind of estimates concerning the Malliavin derivative of the solution.

\begin{proposition}\label{prop:bnd-norm-X-x}
Under the assumptions of Theorem \ref{thm:ex-uniq}, the solution to equation (\ref{eds:alg-form}) satisfies:
\begin{equation}\label{eq:bnd-norm-X-x}
\|X^x\|_{\ga} \le c_{V} \lp 1+ \|\bb^{\1}\|_{\ga}^{1/\ga} +  \|\bb^{\2}\|_{2\ga}^{1/2\ga}\rp,
\end{equation}
where $c_V$ is a constant which only depends on the vector fields $V_0,\ldots,V_n$.
\end{proposition} 

\begin{remark}
This proposition is shown in \cite{FV-bk} by identifying the signature of $B$ with the signature of a certain finite variation process, plus some easy estimates for ordinary differential equations. We have included here a direct elementary proof of (\ref{eq:bnd-norm-X-x}) because we haven't been able to find them in the literature under this form, and mostly because the generalization of our estimates to linear cases will be obvious from the considerations below.
\end{remark}

\begin{proof}[Proof of Proposition \ref{prop:bnd-norm-X-x}]
Theorem \ref{thm:ex-uniq} part (iii) asserts the convergence of the approximation $X^n$ towards $X^x$ as $n\to\infty$. It is thus sufficient to prove relation (\ref{eq:bnd-norm-X-x}) for $X^n$, uniformly in $n$. One can also be easily reduced to prove
\begin{equation*}
\sup_{0\le i<j\le 2^n} \frac{| \der X^n_{t_i^n t_j^n}|}{|t_j^n-t_i^n|^{\ga}} \le 
c_V \lp 1+ \|\bb^{\1}\|_{\ga}^{1/\ga} +  \|\bb^{\2}\|_{2\ga}^{1/2\ga}\rp,
\end{equation*}
which is what we shall proceed to do. We now divide our proof in several steps.

\smallskip

\noindent
\textit{Step 1: Expression for $\der X^n$.}
Set 
\beq\label{eq:def-q-n}
q_{st}^{n}=V_0(X^{n}_{s}) \, (t-s) + V_{i}(X^{n}_{s}) \bb^{\1,i}_{st}
+V_{i_1}V_{i_2}(X^{n}_{s}) \bb^{\2,i_1i_2}_{st},
\eeq
so that $\der X_{t_{i}^{n}t_{j}^{n}}^{n}=\sum_{l=1}^{j-1} q^{n}_{t_{l}^{n}t_{l+1}^{n}}$. For $i<j$, we also construct a dyadic partition $\{\tau_l^k;\, 0\le k\le K, 0\le l\le 2^k\}$ of the set $\{t_i^n,\ldots,t_j^n\}$ inductively in the following way:  set $\tau_0^0=t_i$ and $\tau_1^0=t_j$. Now, if the $\tau_l^k$'s are known, we set $\tau_{2l}^{k+1}=\tau_l^k$. Furthermore, if  $\tau_{2l}^{k+1}=t_m$ and $\tau_{2l+2}^{k+1}=t_{m'}$, then take $\tau_{2l+1}^{k+1}=t_{m^*}$, with $m^*=\lfloor (m+m')/2\rfloor$. This procedure is then non trivial as long as $j-i\ge 2^k$, which means that we stop at $K=\lceil \log_2(j-i)\rceil$.  Here is then a simple example of construction: consider $i=1,j=4$. Then we have:
$$
\tau_0^0=1,\, \tau_1^0=4; \quad
\tau_0^1=1,\, \tau_1^1=2,\, \tau_2^1=4; \quad
\tau_0^2=1,\, \tau_1^2=1,\,\tau_2^2=2,\, \tau_3^2=3,\,\tau_4^2=4.
$$

\smallskip

With these notations in hand, it is easily checked that the relation $\der X_{t_{i}^{n}t_{j}^{n}}^{n}=\sum_{l=1}^{j-1} q^{n}_{t_{l}^{n}t_{l+1}^{n}}$ can also be written as $\der X_{t_{i}^{n}t_{j}^{n}}^{n}=\sum_{k=0}^{2^K-1}q^{n}_{\tau_{k}^{K}\tau_{k+1}^{K}}$. Furthermore,
$$
q_{\tau_{2l}^{k}\tau_{2l+1}^{k}}^{n} + q_{\tau_{2l+1}^{k}\tau_{2l+2}^{k}}^{n}
=q_{\tau_{l}^{k-1}\tau_{l+1}^{k-1}}^n
-(\der q^{n})_{\tau_{2l}^{k}\tau_{2l+1}^{k}\tau_{2l+2}^{k}},
$$
and summing this equality for $k=K$ and $l=0,\ldots,2^{K}-1$ we get
$$
\sum_{l=0}^{2^K-1} q_{\tau_{l}^{K}\tau_{l+1}^{K}}^n
=\sum_{l=0}^{2^K-1} q_{\tau_{l}^{K-1}\tau_{l+1}^{K-1}}^n
-\sum_{l=0}^{2^K-1} (\der q^{n})_{\tau_{2l}^{K-1}\tau_{2l+1}^{K-1}\tau_{2l+2}^{K-1}}.
$$
Iterating, we obtain
\beq\label{eq:exp-delta-y-n}
\der X_{t_{i}^{n}t_{j}^{n}}^{n}=q_{t_{i}^{n}t_{j}^{n}}^{n}
-\sum_{k=1}^{K}\sum_{l=0}^{2^k-1} \der q_{\tau_{2l}^{k}\tau_{2l+1}^{k}\tau_{2l+2}^{k}}^n.
\eeq

\smallskip

\noindent
\textit{Step 2: Expression for $\der q^n$.} Denote by $I$ the identity function on $\R$, so that $\der I_{st}=t-s$. Start then from expression (\ref{eq:exp-delta-y-n}) and use Proposition \ref{difrul} in order to get, for any $s,u,t$ in the dyadic partition,
\begin{equation*}
\der q_{sut}^{n}=
\der V_0(X^n)_{su} \, (t-u)+\der V_i(X^n)_{su} \bb^{1,i}_{ut} 
+ \der\lc V_{i_1}V_{i_2}(X^{n}) \rc_{su} \bb^{\2,i_1i_2}_{ut}
-V_{i_1}V_{i_2}(X^{n}_{s}) \der \bb^{\2,i_1i_2}_{sut},
\end{equation*}
or otherwise stated thanks to convention (\ref{cvpdt}),
\begin{equation}\label{eq:exp-der-q-n}
\der q^{n}=
\der V_0(X^n) \, \der I+\der V_i(X^n) \bb^{1,i}
+ \der\lc V_{i_1}V_{i_2}(X^{n}) \rc \bb^{\2,i_1i_2}
-V_{i_1}V_{i_2}(X^{n}) \der \bb^{\2,i_1i_2}.
\end{equation}
Observe now that one can prove, as for (\ref{eq:dcp-V-i-X}), that 
\begin{equation*}
\der V_i(X^n)= V_{i_1} V_{i}(X^n) \bb^{\1,i_1} +r^{i,n},
\quad\mbox{with}\quad
| r^{i,n}_{st} | \le c_V |\der X_{st}^{n}|^{2}.
\end{equation*}
Thus, according to the fact that $\der\bb^{\2,i_1i_2}=\bb^{\1,i_1}\bb^{\1,i_2}$, one can recast (\ref{eq:exp-der-q-n}) into 
\begin{equation}\label{eq:exp-der-q-n-2}
\der q^{n}=
\der V_0(X^n) \, \der I + r^{i,n} \bb^{1,i} + \der\lc V_{i_1}V_{i_2}(X^{n}) \rc \bb^{\2,i_1i_2}
:= \rho^{1,n} + \rho^{2,n} + \rho^{3,n}.
\end{equation}
Notice that for $j=1,2,3$, the increment $\rho^{j,n}$ lies into $\cac_3^{3\ga}$. Furthermore, it is readily checked that
\begin{equation}\label{eq:bnd-rho-n}
|\rho^{1,n}_{sut}| \le c_V |\der X^{n}_{su}| \, |t-u|^{1+\ga}, \quad
|\rho^{2,n}_{sut}| \le c_V |\der X^{n}_{su}|^{2}  \|\bb^{\1}\|_{\ga}  |t-u|^{\ga}, 
\end{equation}
and $|\rho^{3,n}_{sut}| \le c_V |\der X^{n}_{su}|  \|\bb^{\2}\|_{2\ga}  |t-u|^{2\ga}$.

\smallskip

\noindent
\textit{Step 3: An induction procedure.} Let us consider an integer $\ell\ge 1$ and the quantity:
\begin{equation*}
N_{\ell}^{n}=\sup_{0\le i<j\le \ell} \frac{|\der X^n_{t_i^nt_j^n}|}{|t_{j}^{n}-t_{i}^{n}|^{\ga}}.
\end{equation*}
We localize now our study to an interval of the form $[a,a+\eta]$ with an arbitrary positive number $a$, and $\eta$ small enough. We will prove that if $\eta$ is of order $(1+\|\bb^{\1}\|_{\ga}^{1/2\ga}+\|\bb^{\2}\|_{2\ga}^{1/2\ga})^{-1}$, then $N_{\ell} \le c_V(1+\|\bb^{\1}\|_{\ga}+\|\bb^{\2}\|_{2\ga}^{1/2})$ by induction.

\smallskip

The case $\ell=1$ being trivial, let us assume that the hypothesis is true up to a given $\ell\ge 1$. Take now  $1\le i\le \ell$ and $j=\ell+1$. According to (\ref{eq:exp-delta-y-n}), write 
$$
\der X_{t_{i}^{n}t_{j}^{n}}^{n}-q_{t_{i}^{n}t_{j}^{n}}^{n}=
-\sum_{k=1}^{K}\sum_{l=0}^{2^k-1} \der q_{\tau_{2l}^{k}\tau_{2l+1}^{k}\tau_{2l+2}^{k}}^n.
$$
In the right hand side of this decomposition, all the points $\tau_{2l}^{k},\tau_{2l+1}^{k}$ are of the form $t_p^n$ with $p\le \ell$. Thus (\ref{eq:exp-der-q-n-2}), our bounds (\ref{eq:bnd-rho-n}) on $\rho^{j,n}$ and the induction hypothesis entail
\begin{equation*}
\lln   \der X_{t_{i}^{n}t_{j}^{n}}^{n}-q_{t_{i}^{n}t_{j}^{n}}^{n}\rrn
\le c_V \lp  N_{\ell} + N_{\ell}^{2} \|\bb^{\1}\|_{\ga} + N_{\ell} \|\bb^{\2}\|_{2\ga} \rp |t_{j}^{n}-t_i^n|^{3\ga}.
\end{equation*}
Furthermore, it is obvious from (\ref{eq:def-q-n}) that
\begin{equation*}
\frac{|q_{t_{i}^{n}t_{j}^{n}}^{n}|}{|t_{j}^{n}-t_{i}^{n}|^{\ga}} \le 
c_V \lp 1+  \|\bb^{\1}\|_{\ga} +\|\bb^{\2}\|_{2\ga}|t_{j}^{n}-t_{i}^{n}|^{\ga}\rp.
\end{equation*}
Hence, putting together the last two inequalities, taking into account that we work on an interval of size $\eta$ and that we have chosen $j=\ell+1$, we end up with the following induction relation: $N_{\ell+1}\le F_{\eta}(N_{\ell})$, where the function $F_{\eta}:\R_+\to\R_+$ is defined by
\begin{equation*}
F_{\eta}(\xi)= \xi \vee c_V \lc  \lp 1+  \|\bb^{\1}\|_{\ga} +\|\bb^{\2}\|_{2\ga} \, \eta^{\ga} \rp 
+ \lp 1 +\|\bb^{\2}\|_{2\ga} \, \xi + \|\bb^{\1}\|_{\ga} \xi^2  \rp  \eta^{2\ga}
\rc .
\end{equation*}
By separating the cases $\eta^{\ga}\xi\le 1$ and $\eta^{\ga}\xi > 1$, one can also prove that $F_{\eta}(\xi)\le \xi \vee  \vp_{\eta}(\xi)$, with
\begin{equation*}
\vp_{\eta}(\xi)=  
c_V \lp \|\bb^{\2}\|_{2\ga} \, \eta^{3\ga} + \|\bb^{\1}\|_{\ga} \eta^{2\ga} \rp    \xi^2
+\lp 1+  \|\bb^{\1}\|_{\ga} +\|\bb^{\2}\|_{2\ga} \, \eta^{\ga} \rp 
:= a \, \xi^2 + c.
\end{equation*}
In order to obtain a bound of the form $N_\ell\le M$ which remains valid for all $\ell$, it is now sufficient to have the interval $[0,M]$ left invariant by $ \vp_{\eta}$.

\smallskip

We let the reader check the following elementary fact: whenever $4ac$ is of order 1, the interval $[0,M]$ is left invariant by the application $\xi\mapsto a \, \xi^2 + c$, with $M$ of order $c$. Applying this to $\vp_{\eta}$, we obtain that
\begin{equation}\label{eq:loc-estimate-X-n}
\eta \asymp \lc  1+  \|\bb^{\1}\|_{\ga} +\|\bb^{\2}\|_{2\ga}^{1/2} \rc^{-1/\ga}
\quad \Longrightarrow\quad
N_\ell \lesssim  \|\bb^{\1}\|_{\ga} + \|\bb^{\2}\|_{2\ga}^{1/2}.
\end{equation}

\smallskip

\noindent
\textit{Step 4: Conclusion.}
We have thus obtained that on any interval of length $\eta$ given by (\ref{eq:loc-estimate-X-n}), we have $\|X^n\|_{\ga}\lesssim  \|\bb^{\1}\|_{\ga} + \|\bb^{\2}\|_{2\ga}^{1/2}$. Our claim (\ref{eq:bnd-norm-X-x}) is now easily deduced by dividing an arbitrary interval $[s,t]$ into subintervals of length $\eta$ as in (\ref{eq:global-holder}).

\end{proof}

\subsection{Estimates for the Malliavin derivative}
We are now interested in extending Proposition \ref{Prop: bnd-first-deriv-holder} beyond the Young setting. Recall thus that we are concerned with equation (\ref{eq:ed-def}), which can be written in our algebraic integration setting as
\begin{equation}\label{eq:alg-eq-M}
\der\cm_{st}= \int_s^t \cm_u  \omega_{0} (X^x_u) \, du 
+ \cj_{st}\lp \cm  \, \omega_{i} (X^x) \, dB^{i} \rp, 
\quad 0\le s<t\le T,
\end{equation}
with final condition $\cm_T=V(X_T^x)$. Then we have the following equivalent of  Theorem \ref{thm:ex-uniq}:
\begin{proposition}
Theorem \ref{thm:ex-uniq} holds true for equation (\ref{eq:alg-eq-M}) under Hypothesis \ref{H} and~\ref{hyp:x}. The discretization scheme for $\cm$ can be written as:
\begin{equation}\label{eq:davie-scheme-M}
\der \cm^{n}_{t_k t_{k+1}}= 
\cm^{n}_{t_k} \lc
\frac{\om_0(X^{n}_{t_k})}{2^n} \der I_{t_kt_{k+1}}
+ \om_{i}(X^{n}_{t_k}) \bb^{\1,i}_{t_k t_{k+1}}
+ \lp  \om_{i_1}\om_{i_2}(X^{n}_{t_k}) + V_{i_1}\om_{i_2}(X^{n}_{t_k})  \rp \bb^{\2,i_2 i_1}_{t_k t_{k+1}}
\rc.
\end{equation}
\end{proposition}

\begin{proof}
As in the case of $X^x$, we only justify expression (\ref{eq:davie-scheme-M}). Note that, according to equation (\ref{eq:alg-eq-M}), we have $\der\cm_{st}=\ze_s^l \bb_{st}^{\1,l}+r$, with $\ze_s^l=\cm_s \om_l(X_s^x)$. Hence, invoking Proposition~\ref{difrul} and relation (\ref{eq:dcp-V-i-X}), we get
\begin{equation*}
\der[\cm \, \om_{i_2}]= \cm \lc  \om_{i_1}\om_{i_2}(X^x)+ V_{i_1}\om_{i_2}(X^x)  \rc \bb^{\1,i_{1}}+ r.
\end{equation*}
Expanding now the right hand side of equation (\ref{eq:alg-eq-M}) with the help of Theorem \ref{intg:mdx}, one easily gets
\begin{equation*}
\der \cm= 
\cm \lc
\om_0(X^{x}) \, \der I
+ \om_{i}(X^{x}) \bb^{\1,i}
+ \lp  \om_{i_1}\om_{i_2}(X^{x}) + V_{i_1}\om_{i_2}(X^{x})  \rp \bb^{\2,i_2 i_1}
\rc + r^{\sharp},
\end{equation*}
which gives the desired justification of our scheme (\ref{eq:davie-scheme-M}).

\end{proof}

We are now ready to state and prove our bounds for the process $\cm$:
\begin{proposition}\label{M-bnd-rough}
Let $\cm$ be the unique solution to (\ref{eq:alg-eq-M}) under Hypothesis \ref{H} and~\ref{hyp:x}. Then

\smallskip

\noindent
\emph{(i)} The bound (\ref{eq:bnd-first-deriv-infty}) on $\|\cm\|_{\infty}$ still holds true in our irregular context.

\smallskip

\noindent
\emph{(ii)}  $\|\cm\|_{\ga}$ satisfies inequality (\ref{eq:bnd-norm-X-x}).
\end{proposition}

\begin{proof}
By the continuity property of It\^o's map for $\cm$, the approximation procedure of $\cm$ by $\cm^m$ described in the proof of Theorem \ref{bound MalliavinD} is still valid. The desired bound for $\|\cm\|_{\infty}$ is thus obtained just like in (\ref{eq:bnd-M-antisym}).

\smallskip

Once a bound on $\|\cm\|_{\infty}$ is available, $\|\cm\|_{\ga}$ can be bounded by considering the Davie type scheme (\ref{eq:davie-scheme-M}), along the same lines as for $X^x$. Details are left to the reader for sake of conciseness.

\end{proof}

\subsection{Density upper bound} We finish this section by extending the density estimate and functional inequalities obtained in the smooth case (when $H>\frac{1}{2}$) to the irregular case (when $\frac{1}{3}<H<\frac{1}{2}$).   We first show that in the rough case, we can obtain the smoothness of density of $X_t$ under Hypothesis \ref{H'}.  Indeed, for this purpose, we only need to show the following integrability of  Malliavin matrix.

\begin{theorem}\label{Malliavin-matrix-bnd-rough}
Fix $H\in (\frac{1}{3}, \frac{1}{2})$. Assume Hypothesis \ref{H'}.  Let $\gamma_X$ be the Malliavin matrix of $X^x_1$, we have 
$|\det\gamma_X|^{-1}\in L^{\infty}(\mp)$. Therefore  $X_1^x$ admits a smooth density.

\end{theorem}
\begin{proof}
It suffices to show that there exists $C>0$ such that for all $v\in\mr^d$,
$$v^T\gamma_Xv\geq C|v|^2.$$

Recall that 
$$\gamma_X^{ij}=\langle \mathbf{D}_. X^i, \mathbf{D}_.X^j\rangle_\mathcal{H} $$
and that the $d\times d$ matrix $\cm$ is such that its $j$-th column $\cm^j$ is given by
$$\cm^j_t=\mathrm{D}^j_tX_1.$$
We have
\begin{align*}
v^T\gamma_Xv=&\sum_{j=1}^d\|v^T\cm^j_.\|^2_\mathcal{H}\\
\geq& C\sum_{j=1}^d\int_0^1(v^T\cm^j_s)^2ds
=C\int_0^1 v^T\cm_s\cm^T_sv ds
\geq C\min_{m}\int_0^1 v^T\cm^m_s(\cm^m_s)^Tv ds.
\end{align*}
In the above, we used that $\mathcal{H}\subset L^2([0,1])$ for the second step, and $\cm^m$ is that described in Theorem \ref{bound MalliavinD}.  

On the other hand one can show, according to the description of $\cm^m$ in Theorem \ref{bound MalliavinD}
$$ v^T\cm^m_s(\cm^m_s)^Tv \geq C|v|^2$$
uniformly for some constant $C>0$, if we assume Hypothesis \ref{H'}. The proof is therefore completed.
\end{proof}

Another consequence of the above boundedness of the Malliavin Matrix $\gamma_X$ is the following.

\begin{theorem}
Fix $H\in (\frac{1}{3}, \frac{1}{2})$. Assume Hypothesis \ref{H'}.  Let $p_X(t,y)$  denote the density function of $X_t^x$. There exist constants $c^{(1)}_t$ and $ c^{(2)}_t$ such that
\begin{equation*}
p_X(t,y) \le c_t^{(1)} \exp\lp - c_t^{(2)}   |y|^\delta\rp\quad\quad\ y\in\mr^d,
\end{equation*}
for any $\delta<H$.
\end{theorem}
\begin{proof}
The proof is similar to that of Theorem \ref{thm:main-gauss-bnd}.  We apply Theorem~\ref{M-bnd-rough} and Theorem~\ref{Malliavin-matrix-bnd-rough} to obtain bounds for $\|\mathbf{D}X_t\|_{k,\gamma}^n$ and $\|\det\gamma_{X_t}^{-1}\|_{L^\alpha}^m$.  The tail estimate for $\mp\{X_t>y\}$ is  derived by Proposition \ref{prop:moments-sdes-rough}.  The rest of the proof is then clear.
\end{proof}

\begin{remark}
We are not able to obtain log-Sobolev inequality as in Theorem \ref{log-sobolev} in the rough case, since when $H<\frac{1}{2}$ the Hilbert norm $\mathcal{H}$ is not controlled by $L^\infty([0,1])$.  On the other hand, by reproducing the proof in Theorem \ref{log-sobolev} together with the following interpolation inequality
$$\|\cdot\|_{\mathcal{H}}\leq C(\|\cdot\|_\infty + \|\cdot\|_\gamma),$$
and estimates in Proposition \ref{M-bnd-rough}, we are able to prove the following version of a \emph{log-Sobolev} type inequality: for $T\in[0,1]$ we have
\begin{align*}
\me f(X_T)^2\ln f(X_T)^2- \left(\me f(X_T)^2\right)\ln \left(\me f(X_T)^2\right)\leq C_{p,T}\left(\me|\nabla f(X_T)|^{2p}\right)^{1/p},
\end{align*}
For all $p>1$. Here $C_{p,T}$ is a universal constant independent of $f$.
\end{remark}




\bigskip

 \begin{thebibliography}{99}
 
 \bibitem{BC}
F. Baudoin, L. Coutin:
Operators associated with a stochastic differential equation driven by fractional Brownian motions.
{\it Stoch. Proc. Appl.} {\bf 117} (2007), no. 5,  550--574.

\bibitem{BH} F. Baudoin, M. Hairer: 
A version of H\"ormander's theorem for the fractional Brownian motion. 
{\it Probab. Theory Related Fields} {\bf 139} (2007), no. 3-4, 373--395.

\bibitem{BO} F. Baudoin, C. Ouyang: Small-time kernel expansion for solutions of stochastic differential equations driven by fractional Brownian motions. {\it Stoch. Proc. Appl.} {\bf 121} (2011), no. 4,  759--792.

\bibitem{BN}M. Besal\'{u}, D. Nualart: Estimates for the solution to stochastic differential equations driven by a fractional Brownian motion with Hurst parameter $H\in (\frac{1}{3}, \frac{1}{2})$. {\it Preprint} (2010).


\bibitem{CHL}M. Capitaine, E. Hsu, M. Ledoux: Martingale representation and logarithmic Sobolev inequality. {\it Electronic Communications in Probability}, {\bf 2} (1997), 71--81. 


\bibitem{CF} T. Cass, P. Friz: Densities for rough differential equations under H\"{o}rmander condition. {\it Ann. Math}. To appear. 

\bibitem{CFV} T. Cass, P. Friz, N. Victoir: 
Non-degeneracy of Wiener functionals arising from rough differential equations. 
{\it Trans. Amer. Math. Soc.} {\bf 361} (2009), no. 6, 3359--3371.

\bibitem{CLL} T. Cass, C. Litterer, T. Lyons:
Integrability estimates for Gaussian rough differential equations.
\textit{Arxiv} preprint (2011).

\bibitem{CQ}
L. Coutin, Z. Qian:
Stochastic analysis, rough path analysis and fractional Brownian motions.
{\it Probab. Theory Related Fields}  {\bf 122}  (2002),  no. 1, 108--140.

\bibitem{Da} 
A. Davie:
 Differential equations driven by rough paths: an approach via discrete approximation.  
{\it Appl. Math. Res. Express.}  (2007),  No. 2,  40 pp.

\bibitem{DNT}
A. Deya, A. Neuenkirch, S. Tindel:
A Milstein-type scheme without L\'evy area terms for SDEs driven by fractional Brownian motion.
{\it Ann. Inst. H. Poincaré Probab. Statist.}, to appear.

\bibitem{Dr}
P. Driscoll: Smoothness of density for the area process of fractional Brownian motion. \textit{Arxiv} preprint (2010)

\bibitem{FV-bk}
P. Friz, N. Victoir:
\emph{Multidimensional dimensional processes seen as rough paths.}
Cambridge University Press (2010).

\bibitem{Gu}
M. Gubinelli:
Controlling rough paths.
{\it J. Funct. Anal.} {\bf 216}, 86-140 (2004).

\bibitem{GT} M. Gubinelli, S. Tindel:
Rough evolution equations.
{\it Ann. Probab.}  {\bf 38}  (2010),  no. 1, 1--75.

\bibitem{Ha}
M. Hairer: 
Ergodicity of stochastic differential equations driven by fractional Brownian motion.  
\textit{Ann. Probab.}  \textbf{33}  (2005),  no. 2, 703--758.

\bibitem{HO}
M. Hairer, A. Ohashi: 
Ergodicity theory of SDEs with extrinsic memory.
\textit{Ann. Probab.}  \textbf{35}  (2007),  no. 5, 1950--1977.  


\bibitem{Hsu} E. Hsu: {\it Stochastic Analysis on Manifolds}, Graduate Series in Mathematics, volume 38, Amer. Math. Soc., Providence, RI (2002).

\bibitem{HN}
Y. Hu, D. Nualart:
Differential equations driven by H\"older continuous functions of 
order greater than $1/2$. 
{\it Abel Symp.} {\bf 2}  (2007), 349-413.

\bibitem{KS}
S. Kou, X. Sunney-Xie:
Generalized Langevin equation with
fractional Gaussian noise: subdiffusion within a single protein molecule.
{\it Phys. Rev. Lett.} {\bf 93}, no. 18 (2004).


\bibitem{Ledoux} M. Ledoux: {\it The Concentration of Measure Phenomenon}, Mathematical Surveys and Monographs, volume 89, Amer. Math. Soc., Providence, RI (2001).

\bibitem{LQ}
T. Lyons, Z. Qian:
{\it System control and rough paths.}
Oxford University Press (2002).

\bibitem{MS}
A. Millet, M. Sanz-Solé:
Large deviations for rough paths of the fractional Brownian motion.  
{\it Ann. Inst. H. Poincaré Probab. Statist.}  {\bf 42}  (2006),  no. 2, 245--271. 

\bibitem{NNRT}
A. Neuenkirch, I. Nourdin, A. R\"o\ss ler, S. Tindel :
Trees and asymptotic developments for fractional diffusion processes.
{\it  Ann. Inst. Henri Poincar\'e Probab. Stat.}  {\bf 45} (2009),  no. 1, 157--174.

\bibitem{NTU}
A. Neuenkirch, S. Tindel, J. Unterberger:
Discretizing the L\'evy area. 
{\it Stochastic Process. Appl.}  {\bf 120}  (2010),  no. 2, 223--254.

 \bibitem{Nu06}
D. Nualart: \emph{The Malliavin Calculus and Related
Topics.} Probability and its Applications. Springer-Verlag, 2nd
Edition, (2006).

\bibitem{NR}
D. Nualart,  A. R\v{a}\c{s}canu:
Differential equations driven by fractional Brownian motion.
{\it Collect. Math.} {\bf 53} no. 1 (2002), 55-81.

\bibitem{NS}
D. Nualart, B. Saussereau:
Malliavin calculus for stochastic differential equations driven by a fractional Brownian motion.  {\it Stochastic Process. Appl.}  {\bf 119}  (2009),  no. 2, 391--409. 


\bibitem{SW}
J. Szymanski, M. Weiss:
Elucidating the origin of anomalous diffusion in crowded fluids.
{\it Phys. Rev. Lett.} {\bf 103}, no. 3 (2009).

\bibitem{TBV}
V. Tejedor, O. Bénichou, R. Voituriez, R. Jungmann, F. Simmel, C. Selhuber-Unkel, L. Oddershede and R. Metzle:
Quantitative Analysis of Single Particle Trajectories: Mean Maximal Excursion Method.
\emph{Biophysical J.} {\bf 98}, no. 7 (2010), 1364-1372.

\bibitem{Us}
A.S. \"Ust\"unel: \emph{Analysis on Wiener Space and Applications}. 
Arxiv Preprint (2010).

\bibitem{Za}
M. Z\"ahle:
Integration with respect to fractal functions and stochastic  calculus I.
{\it Probab. Theory Relat. Fields} {\bf 111}  (1998), 333-374.


\end{thebibliography}
 \end{document}